\documentclass[12pt,a4paper]{amsart}

\usepackage{amsmath,amsfonts,amssymb,mathtools,extarrows}
\usepackage{xcolor}
\usepackage{comment}
\usepackage{float}

\usepackage{color}
\usepackage{graphicx}

\usepackage{enumerate}  

\usepackage{cancel}

\usepackage[hmargin=2cm,vmargin=2cm]{geometry}

\usepackage[hypertexnames=false,hyperfootnotes=false,colorlinks=true,linkcolor=blue,%
citecolor=purple,filecolor=magenta,urlcolor=cyan,unicode,linktocpage=true,pagebackref=false]{hyperref}
\usepackage{nameref,zref-xr}     

\usepackage{tikz}

\newcommand{\mbP}{\mathbb P}
\newcommand{\mbZ}{\mathbb Z}
\newcommand{\mbC}{\mathbb C}
\newcommand{\mbR}{\mathbb R}

\newcommand{\oM}{\overline{\mathcal M}}

\newcommand{\og}{\overline g}

\def\cM{{\mathcal{M}}}
\def\oM{{\overline{\mathcal{M}}}}

\def\mbQ{{\mathbb Q}}
\def\d{{\partial}}
\newcommand{\ch}{{\rm ch}}

\newcommand{\eps}{\varepsilon}

\newcommand{\DR}{\mathrm{DR}}

\newcommand{\even}{\mathrm{even}}
\newcommand{\ct}{\mathrm{ct}}

\newcommand{\Coef}{\mathrm{Coef}}

\newcommand{\tv}{\widetilde v}

\DeclareMathOperator{\Aut}{Aut}

\newcommand{\gl}{\mathrm{gl}}

\DeclareMathOperator{\res}{{res}}

\newcommand{\PGL}{\mathrm{PGL}}

\newcommand{\SL}{\mathrm{SL}}

\newcommand{\vir}{\mathrm{vir}}

\newcommand{\Id}{\mathrm{Id}}

\newcommand{\mbE}{\mathbb{E}}

\newcommand{\lb}{\left(}
\newcommand{\rb}{\right)}

\newcommand{\tlambda}{\widetilde{\lambda}}

\DeclareMathOperator{\Hom}{Hom}
\renewcommand{\Re}{\mathrm{Re}}
\newcommand{\Contr}{\mathrm{Contr}}
\newcommand{\tgamma}{\widetilde{\gamma}}
\newcommand{\cc}{\mathrm{c}}

\newtheorem{theorem}{Theorem}[section]

\theoremstyle{definition}

\newtheorem{remark}[theorem]{Remark}
\newtheorem{example}[theorem]{Example}

\usepackage{color}

\numberwithin{equation}{section}

\renewcommand{\gg}[2]{\fill[color=white] (#2) circle(2.5mm) node {\color{black}$\substack{#1}$}; \draw (#2) circle (2.5mm)}
\newcommand{\lab}[4]{\draw (#1)++(#2:#3) node {$\substack{#4}$};}
\newcommand{\leg}[2]{\begin{scope}[shift={(#1)}] \draw (0:0) -- (#2:8mm);\end{scope}}
\newcommand{\lp}[2]{\begin{scope}[shift={(#1)}] \draw (#2:0.3) circle(0.3);\end{scope}}

\usetikzlibrary{decorations.markings}

\newcommand{\legm}[3]{
  \begin{scope}[shift={(#1)}]
    \draw (#2:2.5mm) -- (#2:8mm);

    \node[
      fill=white,
      inner sep=1pt
    ] at (#2:10.5mm) {\color{black}$\substack{#3}$};
  \end{scope}
}

\newcommand{\legmout}[3]{
  \begin{scope}[shift={(#1)}]
    \draw[
      postaction={
        decorate,
        decoration={
          markings,
          mark=at position 0.5 with {\arrow{>}}
        }
      }
    ] (#2:2.5mm) -- (#2:8mm);

    \node[
      fill=white,
      inner sep=1pt
    ] at (#2:10.5mm) {\color{black}$\substack{#3}$};
  \end{scope}
}

\newcommand{\legin}[2]{
  \begin{scope}[shift={(#1)}]
    \draw[
      postaction={
        decorate,
        decoration={
          markings,
          mark=at position 0.6 with {\arrowreversed{>}}
        }
      }
    ] (#2:2.5mm) -- (#2:8mm);
  \end{scope}
}

\newcommand{\legmin}[3]{
  \begin{scope}[shift={(#1)}]
    \draw[
      postaction={
        decorate,
        decoration={
          markings,
          mark=at position 0.5 with {\arrowreversed{>}}
        }
      }
    ] (#2:2.5mm) -- (#2:8mm);

    \node[
      fill=white,
      inner sep=1pt
    ] at (#2:10.5mm) {\color{black}$\substack{#3}$};
  \end{scope}
}

\tikzset{
  midarrow/.style={
    shorten <=2.5mm,
    shorten >=2.5mm,
    postaction={
      decorate,
      decoration={
        markings,
        mark=at position 0.5 with {\arrow{>}}
      }
    }
  }
}

\begin{document}

\title[Integrable systems and moduli spaces of curves]{Cetraro lectures on integrable systems and moduli spaces of curves}

\author{Alexandr Buryak}
\thanks{The work of A.~B. is an output of the research project HSE-BR-2025-059 implemented as a part of the Basic Research Program at the National Research University Higher School of Economics (HSE University). P.~R. is supported by the University of Padova and is affiliated to the INFN under the national project MMNLP and to the INdAM group GNSAGA}

\address{A. Buryak:\newline 
Faculty of Mathematics, National Research University Higher School of Economics, Usacheva str. 6, Moscow, 119048, Russian Federation;\smallskip\newline 
Skolkovo Institute of Science and Technology, Bolshoy Boulevard 30, bld. 1, Moscow, 121205, Russian Federation}
\email{aburyak@hse.ru}

\author{Paolo Rossi}
\address{P. Rossi:\newline Dipartimento di Matematica ``Tullio Levi-Civita'', Universit\`a degli Studi di Padova,\newline
Via Trieste 63, 35121 Padova, Italy}
\email{paolo.rossi@math.unipd.it}

\begin{abstract}
In these lecture notes we review some constructions and results on the relation between topology of the moduli space of stable curves and integrable systems of evolutionary PDEs.
\end{abstract}

\date{\today}

\maketitle

\tableofcontents

\section{Crash course on moduli of curves}

The history of the relation between modern classical integrability (as the theory of systems of evolutionary PDEs possessing an infinite number of symmetries, or a Lax representation, or a (bi)-Hamiltonian structure, etc.) and intersection theory on the moduli space of stable algebraic curves is long enough for the subject to be considered classical. Tracing back its starting moment is difficult, but a turning point which brought the subject much attention from both the phyisics and the mathematics communities is certainly the statement of Witten's conjecture \cite{Wit91} and its proof by Kontsevich \cite{Kon92}.

\medskip

\subsection{Preliminaries}
Let $\oM_{g,n}$ be the moduli space of genus $g$ stable curves with $n$ marked points. In case the reader is not completely familiar with the basic geometric structure of $\oM_{g,n}$, we will attempt a crash course here on the most operative aspects. The truth is that there are excellent sources where one can learn much more for very little extra effort, for instance \cite{Zvo06,CM16, Sch20}.\\

For the time being, let us recall that elements of $\oM_{g,n}$ are closed Riemann surfaces of genus $g$, smooth or with isolated nodal singularities, i.e. singular points with local model $\{xy=0\} \subset \mathbb{C}^2$, and with $n$ distinct and labeled points (the labels will be distinct elements from a set of integers, typically $[n]:=\{1,\ldots,n\}$), up to holomorphic isomorphisms preserving marked points and the set of nodes. Moreover, each curve $C$ representing a point of $\oM_{g,n}$ has to be stable, which means it has a finite number of holomorphic automorphisms (always preserving marked points and the set of nodes). Stability is equivalent to a simple numerical condition on each irreducible component of $C$: suppose an irreducible component $\widetilde{C} \subset C$ has genus $g(\widetilde{C})$ and carries $\nu$ of the original $n$ markings of $C$, $\mu$ nodal intersections with other components of $C$ and $\lambda$ nodal self-intersections, then the condition is $2 g(\widetilde{C})-2 +n(\widetilde{C})>0$, where $n(\widetilde{C})=\nu+\mu+2\lambda$ is referred to as the number of special points of $\widetilde{C}$.

\medskip

$\oM_{g,n}$ carries the structure of a complex orbifold of complex dimension $3g-3+n$. Orbifolds can be seen as generalizations of manifolds where the local model $U\subset \mbC^N$ is replaced by a quotient of $U$ by the action of a finite group $G$ which in general depends on the neighborhood. \cite{Zvo06, Sch20} contain quick and informal, but effective introductions to orbifolds, or their algebraic-geometric counterpart Deligne--Mumford stacks, respectively. The orbifold structure is typical in moduli problems, and in particular the stabilizer of a point $[C] \in \oM_{g,n}$ is isomorphic to the automorphisms group $\Aut(C)$ of a representative $C$.

\medskip

\begin{example}
Genus $0$ smooth closed Riemann surfaces are isomorphic to the complex projective line $\mbP^1$ and have automorphism group the complex projective linear group $\PGL(2,\mbC)$. Given three marked points on $\mbP^1$, there is a unique projective transformation mapping them to, say, $0$, $1$ and $\infty$, after which no residual nontrivial automorphism is left. This means that~$\oM_{0,3}$ is just a point.
\end{example}

\medskip

\subsection{Boundary stratification and stable graphs}
An important feature of $\oM_{g,n}$ is its stratified structure, whose combinatorial description (central for the connection to integrable systems) can be described as follows. First notice that any point $[C]\in \oM_{g,n}$ has an associated dual graph whose vertices correspond to the irreducible components of $C$, whose edges between two vertices (or from a vertex to itself) correspond to nodes at which two irreducible components intersect (or an irreducible component self-intersects), and whose legs (dangling half-edges attached to just one vertex) correspond to marked points. Vertices and legs come marked with their genus and the marked point labeling, respectively, so that we can talk about the total genus $g$ and number of marked legs $n$ of the graph itself, matching the genus and number of marked points of the corresponding curve.

\medskip

Such decorated graphs, dual to stable curves, are called stable graphs. Stable graphs of genus~$g$ with~$n$ marked legs form the set $G_{g,n}$. For a stable graph $\Gamma$ we denote by $V(\Gamma)$, $E(\Gamma)$, $H(\Gamma)$ and $L(\Gamma)$ the sets of its vertices, edges (excluding legs), half-edges (including legs) and legs, respectively. We often use the notations $v_i \in V(\Gamma)$ for the vertex to which the leg $i$ is attached, or $v_h$ for the vertex to which the half-edge $h$ is attached, or $e'$ and $e''$ for the two half-edges forming the edge $e$. The number of loops of a stable graph will be denoted by $h^1(\Gamma)$, so that $g(\Gamma) = \sum_{v \in V(\Gamma)} g(v) + h^1(\Gamma)$. Automorphisms of a stable graph are permutations of its vertices and half-edges which preserve both the way they are attached to each other and the aforementioned genus and legs labeling. Finally, there is a partial order on the set $G_{g,n}$, where $\Gamma\le\Gamma'$ if and only if $\Gamma$ can be obtained from $\Gamma'$ by collapsing some of its edges and modifying the genus labeling in the natural induced way.
\begin{figure}[t]
\begin{center}
\begin{tabular}{ccc}
\raisebox{-0.45\height}{\includegraphics[trim={6cm 0 71.6cm 0},clip,width=0.2\textwidth]{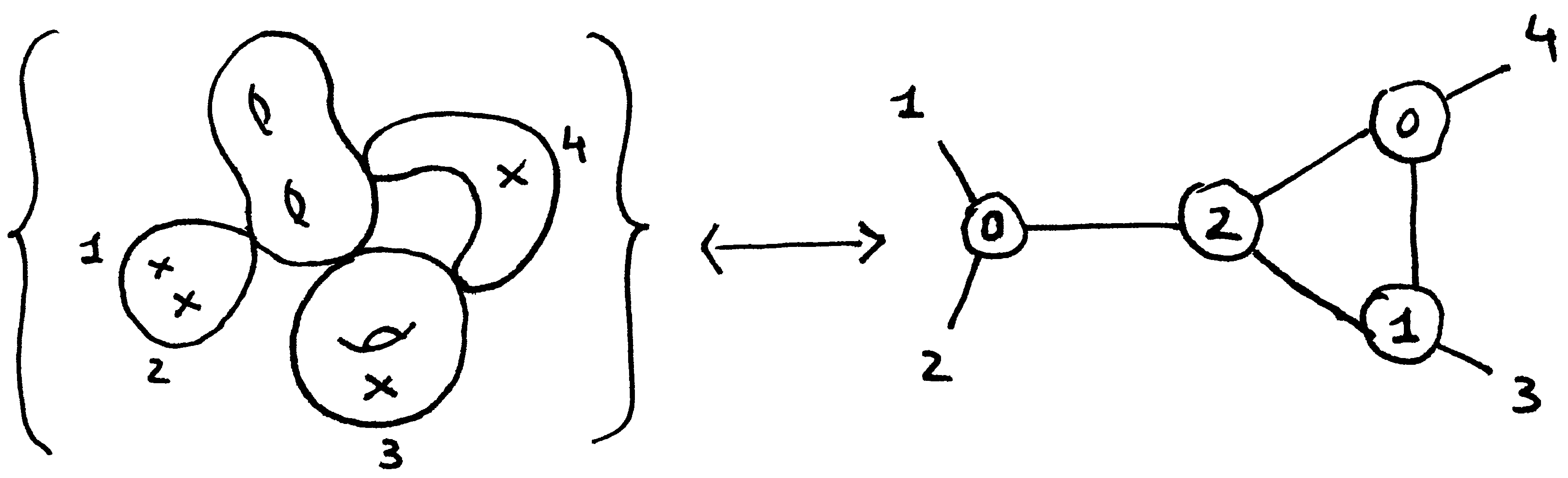}} & \hspace{1cm} &
$\tikz[baseline=-1mm]{\coordinate (A) at (0,0);\coordinate (B) at (12mm,0);\coordinate (C) at (22mm,-5mm);\coordinate (D) at (22mm,5mm);\draw (A)--(B)--(C);\draw (B)--(D)--(C);\legm{A}{120}{1};\legm{A}{-120}{2};\legm{C}{0}{3};\legm{D}{0}{4};\gg{0}{A};\gg{2}{B};\gg{1}{C};\gg{0}{D};}$
\end{tabular}
\end{center}
\caption{A stable curve in $\oM_{4,4}$ and the corresponding stable graph.}
\end{figure}

\medskip

To each stable graph $\Gamma$ we can, vice versa, associate the locus $\cM^{\Gamma} \subset \oM_{g,n}$ of equivalence classes of stable curves $[C]$ whose dual graph is $\Gamma$ and its closure $\oM^{\Gamma}$, which is a closed suborbifold of complex codimension equal to the number of edges of $\Gamma$ (nodes of the curves in $\cM^{\Gamma}$). For any stable graph $\Gamma$, the locus $\oM^{\Gamma} \setminus \cM^{\Gamma}$ is the union of all loci $\cM^{\Gamma'}$, where $\Gamma'\in G_{g,n}$ is strictly bigger than $\Gamma$ with respect to the partial order introduced above. Dually, this means that compactifying $\cM^{\Gamma}$ to $\oM^{\Gamma}$ consists in allowing curves with futher nodal singularities than the ones originally encoded by the edges of $\Gamma$. In particular, the graph with a single genus $g$ vertex and no edges corresponds to $\oM_{g,n}$ itself and the graphs with just one edge correspond to the irreducible components of a divisor with normal crossings representing the boundary of~$\oM_{g,n}$ (i.e. the entire locus of non-smooth stable curves). Given two stable graphs $\Gamma_1,\Gamma_2\in G_{g,n}$, the intersection $\oM^{\Gamma_1} \cap \oM^{\Gamma_2}$ is the union of all loci $\oM^{\Gamma}$ such that $\Gamma\ge\Gamma_1$ and $\Gamma\ge\Gamma_2$.

\medskip

\begin{example}\label{example:M11}
Genus $1$ smooth closed Riemann surfaces (called elliptic curves) are classified by their elliptic modulus: each curve is isomorphic to a torus $E_\tau = \mbC / \Lambda_\tau$ where $\Lambda_\tau=\{m+n\tau | m,n \in \mbZ\}$, for some $\tau$ in the upper half complex plane $\mathbb{H}$. A complex torus is a group, so its automorphism group contains the torus itself as the group of translations, but marking one point kills these translations entirely and all extra symmetries form a finite group. A smooth elliptic curve with a marked point $(E,x_1)$ is isomorphic to $(E_\tau,0)$, and two tori $(E_{\tau_1},0)$ and $(E_{\tau_2}$,0) are isomorphic if and only if $\tau_2 = \frac{a\tau_1 + b}{c\tau_1 + d}$ for some $\begin{pmatrix}a & b \\ c & d \end{pmatrix} \in \SL(2,\mbZ)$. Moreover~$\Aut(E_\tau,0)$ is isomorphic to the stabilizer of $\tau \in \mathbb{H}$ with respect to this $\SL(2,\mbZ)$ action. In other words $\cM_{1,1}\simeq \mathbb{H}/\SL(2,\mbZ)$ as complex orbifolds. The fundamental domain of the action of $\SL(2,\mbZ)$ is given by the shaded region in the familiar picture on Fig.~\ref{fig:fundamental domain}.
\begin{figure}[t]
\begin{center}
\begin{tikzpicture}[scale=2]
    \fill[blue!10] (0.5, 2.5) -- (0.5, {sqrt(3)/2}) arc (60:120:1) -- (-0.5, 2.5) -- cycle;

    \draw[->] (-1.5, 0) -- (1.5, 0) node[right] {$\text{Re}(\tau)$};
    \draw[->] (0, 0) -- (0, 2.8) node[above] {$\text{Im}(\tau)$};

    \draw[dashed] (1,0) arc (0:180:1);
    
    \draw[dashed] (0.5, 0) -- (0.5, {sqrt(3)/2});
    \draw[dashed] (-0.5, 0) -- (-0.5, {sqrt(3)/2});

    \draw[thick] (0.5, {sqrt(3)/2}) -- (0.5, 2.5);
    \draw[thick] (-0.5, {sqrt(3)/2}) -- (-0.5, 2.5);
    \draw[thick] (0.5, {sqrt(3)/2}) arc (60:120:1);

    \draw (-0.5, 0.05) -- (-0.5, -0.05) node[below] {$-\frac{1}{2}$};
    \draw (0.5, 0.05) -- (0.5, -0.05) node[below] {$\frac{1}{2}$};
    \draw (1, 0.05) -- (1, -0.05) node[below] {$1$};
    \draw (-1, 0.05) -- (-1, -0.05) node[below] {$-1$};

    \fill (0, 1) circle (0.7pt) node[below right] {$i$};
    \fill (-0.5, {sqrt(3)/2}) circle (0.7pt) node[above left] {$\rho$};
    \fill (0.5, {sqrt(3)/2}) circle (0.7pt) node[above right] {$\rho+1$};
\end{tikzpicture}
\end{center}
\caption{The fundamental domain of the action of $\SL(2,\mbZ)$ on $\mathbb{H}$.}
\label{fig:fundamental domain}
\end{figure}
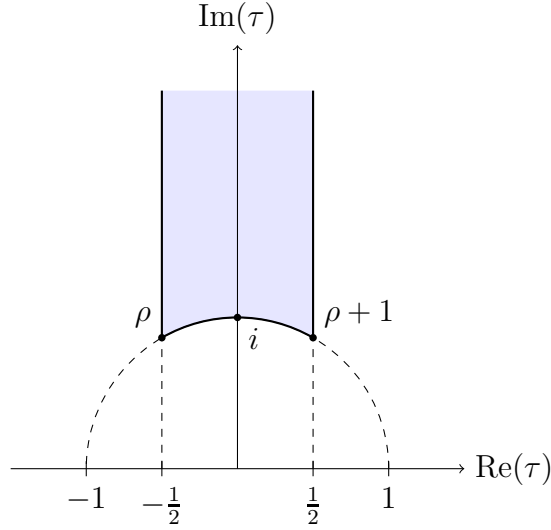
$\cM_{1,1}$ is then represented by this shaded region, where the two vertical lines at $\Re(\tau)=\pm \frac{1}{2}$ are topologically identified and so are the two arcs between the points $\rho$ and $i$ and between the points $\rho+1$ and $i$. Thus $\cM_{1,1}$ is topologically a disk, and we can rely on the above combinatorial description of boundary strata to understand $\oM_{1,1}$. It's easy to see that there are only two stable graphs with $(g,n)=(1,1)$, namely,
$$\tikz[baseline=-1mm]{\coordinate (A) at (0,0);\legm{A}{180}{1};\gg{1}{A};} \quad \text{and}
\quad
\tikz[baseline=-1mm]{\coordinate (A) at (0,0);\lp{A}{0};\legm{A}{180}{1};\gg{0}{A};}\,.$$
The first graph corresponds to the main stratum $\cM_{1,1}$, and the second graph corresponds to the boundary divisor consisting of a single point (representing the only elliptic curve with one marking and one nodal self-intersection), thus compactifying $\cM_{1,1}$ to $\oM_{1,1}$, which is topologically a sphere. The orbifold structure of $\oM_{1,1}$ is encoded in the $\SL(2,\mbZ)$-stabilizers: one immediately checks that for the generic point this is $\mbZ_2=\{\pm \Id\}$ (all elliptic curves have an nontrivial elliptic involution), for $\tau=i$ it is $\mbZ_4$ (the symmetry group of the square lattice) and for $\tau=\rho$ it is $\mbZ_6$ (the symmetry group of the hexagonal lattice).
\end{example}

\medskip

For a stable graph $\Gamma$, let $\oM_\Gamma:=\prod_{v \in V(\Gamma)} \oM_{g(v),n(v)}$ and consider the map $\xi_\Gamma: \oM_\Gamma \to \oM_{g,n}$ gluing stable curves at marked points as prescribed by the edges of the graph $\Gamma$. The image of $\oM_\Gamma$ through $\xi_\Gamma$ is exactly $\oM^{\Gamma}$ but the degree of this map is not always one, in fact, for the fundamental classes, we have $[\oM^\Gamma] = \frac{1}{|\Aut(\Gamma)|} \xi_{\Gamma *} [\oM_\Gamma]$. In formulae, when we represent a stratum by a picture of the corresponding stable graph $\Gamma$, we mean $\xi_{\Gamma *} [\oM_\Gamma]$, leaving out the factor $\frac{1}{|\Aut(\Gamma)|}$.

\medskip

\begin{remark}
There is one thing we should point out about cohomology, intersection theory and Poincaré duality on an orbifold $X$ of complex dimension $\dim X$. We consider singular homology and cohomology with $\mbQ$ coefficients, which are only sensitive to the topological structure of the space. Singular cohomology is isomorphic to the orbifold version of de Rham cohoomlogy, where, in particular a differential form $\alpha$ on $X$ is locally described as a $G$ invariant differential form $\widehat{\alpha}$ on a each orbifold chart $\pi:U\to U/G$. The integral $\int_\gamma \alpha$ over a sufficiently regular chain $\gamma\subset X$ has local contribution $\frac{1}{|G|} \int_{\pi^{-1}(\gamma)}\widehat{\alpha}$. For a compact orbifold $X$, the fundamental class $[D]$ of a sufficiently regular $k$-cycle $D$ (compact sub-orbifold) can then be defined, as usual, via integration over $D$, as an element of $\Hom(H^{k}(X),\mbQ) \simeq H_{k}(X)$. The Poincaré pairing $H^k(X) \otimes H^j(X)\to\mbQ$ given by $(\omega\otimes\eta)\mapsto\int_X \omega \wedge \eta$ is perfect and Poincaré duality is the isomorphism $H^k(X) \to H_{2\dim X - k}(X): \alpha \mapsto [X]\frown \alpha$, so that, in particular it's still true that $1\mapsto [X]$ under Poincaré duality. Let's agree to suppress the notational distinction between Poincaré dual homology and cohomology classes from here on.
\end{remark}

\begin{example}
Let $\Gamma = \tikz[baseline=-1mm]{\coordinate (A) at (0,0);\lp{A}{0};\legm{A}{180}{1};\gg{0}{A};}$ be the second graph appearing in Example \ref{example:M11}. Then $|\Aut(\Gamma)|=2$ (the only nontrivial automorphism exchanges the half-edges of the loop), $\xi_{\Gamma}\colon\oM_{0,3} \to \oM_{1,1}$ with $[\oM_{0,3}]=[\text{pt}]=1 \in H^0(\oM_{0,3})$, and $1=[\oM^\Gamma]= [\text{pt}] =\frac{1}{2} (\xi_\Gamma)_*(1) \in H^2(\oM_{1,1})$ (where~$[\text{pt}]$ denotes the class of the generic point of a space).
\end{example}

\begin{remark}
Natural partial compactifications of $\cM_{g,n}$, contained in $\oM_{g,n}$, are studied in various contexts. In these notes we will be interested in considering in particular the moduli space $\cM_{g,n}^\ct \subset \oM_{g,n}$ of stable curves of compact type, which is the union of strata of $\oM_{g,n}$ whose dual graph is a tree, i.e. only smooth and singular curves with separating nodes are admitted. This space is not compact and, in order to keep using Poicar\'e duality, we will use cohomology with compact support (or Borel--Moore homology). Since the gluing maps $\xi_\Gamma$ for stable trees $\Gamma$ are proper we still have well defined pushforward maps $(\xi_\Gamma)_*$ in ordinary cohomology.
\end{remark}

\medskip

\subsection{Universal curve}
Consider the map $\pi\colon\cM_{g,n+1} \to \cM_{g,n}$ which forgets the $(n+1)$-st marked point. Some care is neeeded when extending the map $\pi$ to $\oM_{g,n+1}$: when an irreducible component of a curve $C$ from $\oM_{g,n+1}$ is rendered unstable upon forgetting $x_{n+1}$, $\pi$ contracts the whole component to a point. By the numerical condition for stability we have seen above, this can only happen if the irreducible component has genus $0$ and either one marked point and two nodes or two marked points and one node. In the first case the two nodes collapse into one node upon forgetting the marking, while in the second case the surviving marked point collapses to the position of the node. The fibration $\pi$ is referred to as the universal curve over $\oM_{g,n}$ since its fiber over $[C]$, up to the orbifold structure, is $C/\Aut(C)$ with its $n$ original markings, seen as the set of points at which the extra $(n+1)$-st marking can be added back up to isomorphism (an automorphism $\phi$ of $(C,x_1,\ldots,x_n)$ is an isomorphism between $(C,x_1,\ldots,x_n,x_{n+1})$ and $(C,x_1,\ldots,x_n,\phi(x_{n+1}))$).

\medskip

Notice that a smooth unmarked point $x$ of the fiber $\pi^{-1}([C])$, as an element of $\oM_{g,n+1}$, corresponds to a curve isomorphic to $C$ with an extra marking $x_{n+1}$ at the position $x \in C$. A marked point $x_i\in \pi^{-1}([C])$ corresponds to the curve obtained from $C$ by replacing $x_i$ with a new node $\nu$ to which an extra genus $0$ component carrying the markings $x_i$ and $x_{n+1}$ and the node $\nu$ is attached, and a node of $\nu \in \pi^{-1}([C])$ corresponds to the curve obtained from $C$ by replacing $\nu$ with two nodes $\nu_1$ and $\nu_2$, one on each branch of $C$ meeting at $\nu$, to which an extra genus $0$ component carrying the markin $x_{n+1}$ and the nodes $\nu_1$ and $\nu_2$ is attached. These extra genus $0$ components are exactly the ones that get contracted upon forgetting $x_{n+1}$.

\medskip

The notation $\oM_{g,n+1} = \overline{\mathcal{C}}_{g,n}$ and $\pi\colon\overline{\mathcal{C}}_{g,n}\to \oM_{g,n}$ is used to stress the interpretation of the forgetful morpshism as universal curve. $\pi$ has $n$ natural sections $\sigma_i$, $1 \leq i \leq n$, given by the the $n$ marked points on each fiber, i.e. $\sigma_i([C]) = x_i \in \pi^{-1}([C])$. Since markings of stable curves are always distinct by definition, the images $D_{i,n+1}=\sigma_i(\oM_{g,n})$ do not intersect. The above discussion of special points of the fibers of $\pi$ implies in particular that the locus $D_{i,n+1}$ is the component of the boundary divisor of $\oM_{g,n+1}$ given by
$$D_{i,n+1} := \sigma_i(\oM_{g,n})= \tikz[baseline=-1mm]{\coordinate (A) at (0,0);\coordinate (B) at (12mm,0);\draw (A)--(B);\legm{A}{120}{i};\legm{A}{-120}{n+1};\leg{B}{-60};\leg{B}{-120};\gg{0}{A};\gg{g}{B};\lab{B}{-90}{6.9mm}{...};\lab{B}{-90}{10mm}{\underbrace{\phantom{aaa}}_{\text{$[n]\setminus \{i\}$}}};}$$
More in general, for $1\leq i,j \leq n+1$ with $i\neq j$, we denote
$$D_{i,j} := \tikz[baseline=-1mm]{\coordinate (A) at (0,0);\coordinate (B) at (12mm,0);\draw (A)--(B);\legm{A}{120}{i};\legm{A}{-120}{j};\leg{B}{-60};\leg{B}{-120};\gg{0}{A};\gg{g}{B};\lab{B}{-90}{6.9mm}{...};\lab{B}{-90}{10mm}{\underbrace{\phantom{aaa}}_{\text{$[n+1]\setminus \{i,j\}$}}};}$$
Let $\Delta$ be the set formed by the union of nodal points of each fiber of $\pi$. We can give an explicit local model for $\pi$ around a point $\nu \in \Delta$ as follows. There exist local coordinates $(T_1, \ldots, T_{3g-4+n}, x, y)$ in a neighborhood of $\nu$ and $(t_1, \ldots, t_{3g-4+n}, t)$ in a neighborhood of $\pi(\nu)$ such that  $t_i \circ \pi=T_i$, $t\circ \pi = xy$ and $x(\nu)=y(\nu)=0$ (see \cite{Zvo06}).

\medskip

\begin{example}
To construct $\oM_{0,4}$ we can realize it as the universal curve over $\oM_{0,3}=\text{pt}$. There is just one fiber isomorphic to $\mbP^1$ with three markings (the sections $D_{i,4}$, $1 \leq i \leq 3$) which correspond to the boundary divisors of $\oM_{0,4}$. Notice how these boundary divisors all represent the same class $[\text{pt}] \in H^2(\oM_{0,4}) = H^2(\mbP^1)$, i.e
\begin{equation}\label{eq:WDVV}
\tikz[baseline=-1mm]{\coordinate (A) at (0,0);\coordinate (B) at (12mm,0);\draw (A)--(B);\legm{A}{120}{1};\legm{A}{-120}{2};\legm{B}{60}{3};\legm{B}{-60}{4};\gg{0}{A};\gg{0}{B};} =
\tikz[baseline=-1mm]{\coordinate (A) at (0,0);\coordinate (B) at (12mm,0);\draw (A)--(B);\legm{A}{120}{1};\legm{A}{-120}{3};\legm{B}{60}{2};\legm{B}{-60}{4};\gg{0}{A};\gg{0}{B};} =
\tikz[baseline=-1mm]{\coordinate (A) at (0,0);\coordinate (B) at (12mm,0);\draw (A)--(B);\legm{A}{120}{1};\legm{A}{-120}{4};\legm{B}{60}{3};\legm{B}{-60}{2};\gg{0}{A};\gg{0}{B};}
\in H^2(\oM_{0,4})
\end{equation}
This simple cohomological equality is central in the relation between intersection theory of $\oM_{g,n}$ and integrability and is called WDVV (Witten--Dijkgraaf--Verlinde--Verlinde) relation because of its role in topological field theory.
\end{example}

\medskip

\begin{example}
One can proceed in a similar way to describe $\oM_{0,5}$ as the universal curve over~$\oM_{0,4}$. Over the complement of $\bigcup_{i=1}^3 D_{i,4}\subset \oM_{0,4}$ we get a $\mbP^1$-bundle, but the fibers over~$D_{i,4}$ are isomorphic to a nodal curve of genus~$0$, with dual graph given by each term in equation \eqref{eq:WDVV}. These singular fibers are the irreducible components $D_{i,4} \subset \oM_{0,5}$, $1 \leq i \leq 3$, of the boundary divisor, but there are more irreducible components, namely the $D_{i,j}$, with $1 \leq i<j \leq 5$. We can list the corresponding stable graphs:
\begin{figure}[H]
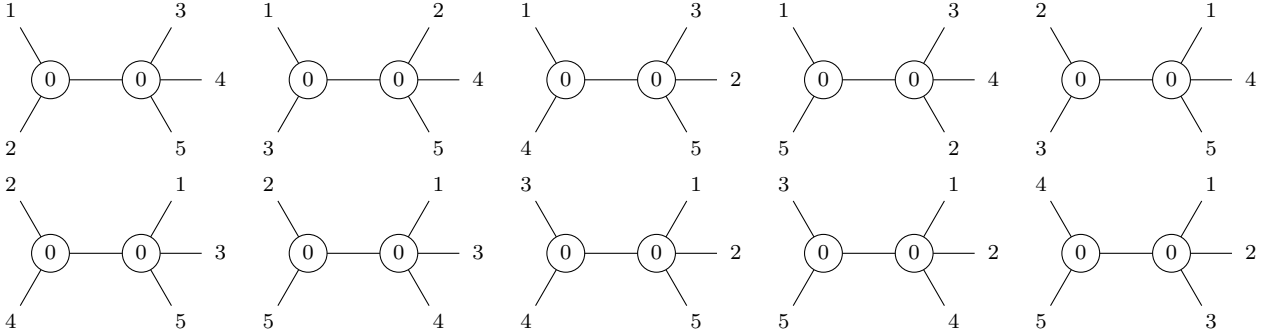

$$
\tikz[baseline=-1mm]{\coordinate (A) at (0,0);\coordinate (B) at (12mm,0);\draw (A)--(B);\legm{A}{120}{1};\legm{A}{-120}{2};\legm{B}{60}{3};\legm{B}{0}{4};\legm{B}{-60}{5};\gg{0}{A};\gg{0}{B};} \quad
\tikz[baseline=-1mm]{\coordinate (A) at (0,0);\coordinate (B) at (12mm,0);\draw (A)--(B);\legm{A}{120}{1};\legm{A}{-120}{3};\legm{B}{60}{2};\legm{B}{0}{4};\legm{B}{-60}{5};\gg{0}{A};\gg{0}{B};} \quad
\tikz[baseline=-1mm]{\coordinate (A) at (0,0);\coordinate (B) at (12mm,0);\draw (A)--(B);\legm{A}{120}{1};\legm{A}{-120}{4};\legm{B}{60}{3};\legm{B}{0}{2};\legm{B}{-60}{5};\gg{0}{A};\gg{0}{B};} \quad
\tikz[baseline=-1mm]{\coordinate (A) at (0,0);\coordinate (B) at (12mm,0);\draw (A)--(B);\legm{A}{120}{1};\legm{A}{-120}{5};\legm{B}{60}{3};\legm{B}{0}{4};\legm{B}{-60}{2};\gg{0}{A};\gg{0}{B};} \quad
\tikz[baseline=-1mm]{\coordinate (A) at (0,0);\coordinate (B) at (12mm,0);\draw (A)--(B);\legm{A}{120}{2};\legm{A}{-120}{3};\legm{B}{60}{1};\legm{B}{0}{4};\legm{B}{-60}{5};\gg{0}{A};\gg{0}{B};}
$$
$$
\tikz[baseline=-1mm]{\coordinate (A) at (0,0);\coordinate (B) at (12mm,0);\draw (A)--(B);\legm{A}{120}{2};\legm{A}{-120}{4};\legm{B}{60}{1};\legm{B}{0}{3};\legm{B}{-60}{5};\gg{0}{A};\gg{0}{B};} \quad
\tikz[baseline=-1mm]{\coordinate (A) at (0,0);\coordinate (B) at (12mm,0);\draw (A)--(B);\legm{A}{120}{2};\legm{A}{-120}{5};\legm{B}{60}{1};\legm{B}{0}{3};\legm{B}{-60}{4};\gg{0}{A};\gg{0}{B};} \quad
\tikz[baseline=-1mm]{\coordinate (A) at (0,0);\coordinate (B) at (12mm,0);\draw (A)--(B);\legm{A}{120}{3};\legm{A}{-120}{4};\legm{B}{60}{1};\legm{B}{0}{2};\legm{B}{-60}{5};\gg{0}{A};\gg{0}{B};} \quad
\tikz[baseline=-1mm]{\coordinate (A) at (0,0);\coordinate (B) at (12mm,0);\draw (A)--(B);\legm{A}{120}{3};\legm{A}{-120}{5};\legm{B}{60}{1};\legm{B}{0}{2};\legm{B}{-60}{4};\gg{0}{A};\gg{0}{B};} \quad
\tikz[baseline=-1mm]{\coordinate (A) at (0,0);\coordinate (B) at (12mm,0);\draw (A)--(B);\legm{A}{120}{4};\legm{A}{-120}{5};\legm{B}{60}{1};\legm{B}{0}{2};\legm{B}{-60}{3};\gg{0}{A};\gg{0}{B};}
$$
\caption{The irreducible components of the boundary divisor of $\oM_{0,5}$.}
\end{figure}
We know how these components intersect each other from the combinatorics of collapsing edges in the corresponding stable graphs. All intersection points are codimension $2$ strata in the boundary, so they correspond to a tree with two edges. Collapsing one or the other edge we obtain each of the two intersecting irreducible divisors.

Geometrically $\oM_{0,5}$ is isomorphic to $(\mbP^1)^2$ blown up at $(0,0)$, $(1,1)$ and $(\infty,\infty)$. The two projections to each of the two copies of $\mbP^1$ are the maps $\pi_4,\pi_5\colon\oM_{0,5}\to\oM_{0,4}\simeq \mbP^1$ forgetting the $4$-th or the $5$-th marked point, respectively. If we take $0,1,\infty \in \mbP^1$ to represent the boundary divisors $D_{2,3}$, $D_{1,3}$, $D_{1,2} \subset \oM_{0,4}$, respectively, then the exceptional divisors represent $D_{2,3}$, $D_{1,3}$, $D_{1,2} \subset \oM_{0,5}$. The fibers of $\pi_4$ and $\pi_5$ over $0$, $1$ and $\infty$ are nodal with two rational components, as they include the exceptional divisors. Explictly, let $k=4,5$, then $\pi_k^{-1}(D_{2,3})=D_{2,3} + D_{1,9-k}$, $\pi_k^{-1}(D_{1,3})=D_{1,3} + D_{2,9-k}$, $\pi_k^{-1}(D_{1,2})=D_{1,2} + D_{3,9-k}$. The rational components of the singular fibers of $\pi_k$ that are not exceptional divisors are the sections $D_{1,9-k}$, $D_{2,9-k}$, $D_{3,9-k}$ of $\pi_{9-k}$. The preimage in the blowup $\oM_{0,5}$ of the diagonal of $(\mbP^1)^2$ minus the three exceptional divisors, provides us with a fourth section of both $\pi_4$ and $\pi_5$, namely $D_{4,5}$.
\end{example}

\medskip

\subsection{Tautological bundles}
Now consider the relative cotangent bundle $\omega_{g,n} \to \overline{\mathcal{C}}_{g,n}\setminus \Delta$ with respect to $\pi$, i.e. the dual bundle to the tangent bundle to the smooth (non-nodal) points of the $\pi$-fiber foliation. To extend $\omega_{g,n}$ over $\Delta$ to the whole of $\oM_{g,n}$ we need to use the local model of $\pi$ around a fiber nodal point $\nu\in \Delta$. Then $\omega_{g,n}$ is defined in a neighborhood of $\nu$ as the line bundle generated by the meromorphic differentials $\frac{dx}{x}$ and $\frac{dy}{y}$ modulo $\langle \frac{d(xy)}{xy} = \frac{dx}{x}+\frac{dy}{y} \rangle$. This coincide with the relative cotangent bundle's definition in the complement of $\nu$, as the line bundle generated by $dT_1\ldots,dT_{3g-4+n},dx,dy$ modulo $\langle \pi^*(dt_1)=dT_1,\ldots,\pi^*(dt_{3g-4+n})=dt_{3g-4+n},\pi^*(dt)=d(xy)\rangle$.

\medskip

For each point $[C] \in \oM_{g,n}$, consider the cotangent complex line $T_{x_i}^*C$ to $C$ at the $i$-th marked point, where $1 \leq i \leq n$. Varying $[C]$ in $\oM_{g,n}$, these lines form a line bundle which we denote by $L_i$ and its first Chern class is denoted by $\psi_i=c_1(L_i) \in H^2(\oM_{g,n})$. A more precise definition is simply $L_i:=\sigma_i^* \omega_{g,n}$. By construction the restriction of $L_i$ to $\oM^\Gamma$ is only sensitive to the geometry of the irreducible component $C_{v_i}$ of $C$ carrying the marked point~$x_i$, where $v_i \in V(\Gamma)$, which means that $\xi_\Gamma^*L_i = p_{v_i}^*L_i$, where $p_{v_i}$ is the cartesian projection $p_{v_i}\colon\prod_{v \in V(\Gamma)} \oM_{g(v),n(v)} \to \oM_{g(v_i),n(v_i)}$. In particular $\xi_{\Gamma}^* \psi_i = p_{v_i}^* \psi_i$.

\medskip

Since the cotangent line at $x_i$ is not influenced by the presence of an extra marked point $x_{n+1}$ at another position on $C$, we have that, $\pi^*L_i = L_i$ as an equality of bundles over $\cM_{g,n+1}$. When trying to extend this equality to $\oM_{g,n}$ we notice that contractions of irreducible components rendered unstable upon forgetting $x_{n+1}$ can instead influence the local geometry around $x_i$ and this happens exactly when $x_{i}$ and $x_{n+1}$ lie on the same contracted component. This is why we have $ L_i = \pi^*L_i + D_{i,n+1}$ as an equality of bundles over $\oM_{g,n+1}$ and, correspondingly, $\psi_i = \pi^*\psi_i + D_{i,n+1} \in H^2(\oM_{g,n+1})$, for their Chern classes (here, and in what follows, we use Poincaré duality in $\oM_{g,n+1}$ to identify homology cycles with cohomology classes and vice versa).

\medskip

\begin{remark}\label{remark:psi-pullback}
Repeated application of the formula $\psi_i = \pi^*\psi_i + D_{i,n+1}$ can be used to reduce~$\psi_i$ on $\oM_{g,n}$ to an explicit expression depending only on the pullback of $\psi_i$ on $\oM_{g,1}$, where the only marked point has the label $i$, unless $g=0$ in which case we have to stop at $\oM_{0,3}$, because of stability, and where the three marking are labeled, say, $i,j,k \in [n]$. But, since $\oM_{0,3}$ is just a point, $\psi_i=0$ on $\oM_{0,3}$ and we obtain
$$\psi_i = \sum_{\substack{I\sqcup J = [n] \\ i\in I,\, j,k\in J}} \tikz[baseline=-1mm]{\coordinate (A) at (0,0);\coordinate (B) at (12mm,0);\draw (A)--(B);\leg{A}{-60};\leg{A}{-120};\leg{B}{-60};\leg{B}{-120};\gg{0}{A};\gg{0}{B};\lab{A}{-90}{6.9mm}{...};\lab{A}{-90}{10mm}{\underbrace{\phantom{aaa}}_{\text{$I$}}};\lab{B}{-90}{6.9mm}{...};\lab{B}{-90}{10mm}{\underbrace{\phantom{aaa}}_{\text{$J$}}};} \quad \in H^2(\oM_{0,n})$$
\end{remark}

\medskip

There are other interesting equations involving the divisors $D_{i,n+1}$ and $\psi$ classes. For instance, we have already seen that $D_{i,n+1}$ and $D_{j,n+1}$ do not intersect, or, in the cohomology ring, $D_{i,n+1}D_{j,n+1}=0 \in H^4(\oM_{g,n+1})$.  Moreover, since $\oM_{0,3}$ is just a point, $L_i|_{D_{i,n+1}}$ is the trivial line bundle, i.e. $\psi_i D_{i,n+1} = 0$. Since $\pi_*D_{j,n+1} =[ \oM_{g,n}]$ (as $D_{j,n+1}$ is a section of $\pi$) we get $\pi_*(D_{j,n+1}\pi^*\psi_i) = \psi_i \in \oM_{g,n}$. Multiplying the equation $\psi_i = \pi^*\psi_i + D_{i,n+1}$ by $D_{i,n+1}$ we thus obtain $D_{i,n+1}^2 = - D_{i,n+1} \pi^*\psi_i$ and more in general $D_{i,n+1}^k=(-1)^{k-1} D_{i,n+1} \pi^*\psi_i^{k-1}$ for all $k\geq 1$, so that $\psi_i^{k-1} = (-1)^{k-1} \pi_*(D_{i,n+1}^k)$. Since $\pi^* \psi_i^k=(\psi_i - D_{i,n+1})^k=$ $\psi_i^k + (-1)^k D_{i,n+1}^k$ $=\psi_i^k - \pi^* \psi_i^{k-1} D_{i,n+1}$ we obtain by pushforward that $\pi_* \psi_i^k = \psi_i^{k-1}$. More in general, but with similar arguments,
$$\pi^*(\psi_1^{k_1} \ldots\psi_n^{k_n}) = \psi_1^{k_1} \ldots\psi_n^{k_n} - \sum_{\substack{1\le i\le n\\k_i\ge 1}} \pi^*\big(\psi_i^{k_i-1} \prod_{j \neq i} \psi_j^{k_j}\big)D_{i,n+1}.$$
From this last formula we easily deduce, in particular, that $\pi_*(\psi_1^{k_1} \ldots\psi_n^{k_n}) =\sum_{\substack{1\le i\le n\\k_i\ge 1}} \psi_i^{k_i-1} \prod_{j \neq i} \psi_j^{k_j}$.\\

\medskip

Since we have been considering the effect of pushing forward powers of psi classes at a marked point with respect to morphisms forgetting another point, it is perhaps the right moment to introduce kappa classes as $\kappa_a = \pi_* \psi_{n+1}^{a+1}$, for integers $a\geq 0$. Since they are the fruit of pushing forward powers of psi classes at a point with respect to the morphism forgetting that point, kappa classes are not attached to marked points like psi classes are. In particular, when restricted to boundary strata, they behave according to the formula $\xi_\Gamma^* \kappa_a = \sum_{v \in V(\Gamma)} p_{v}^* \kappa_a$.

\medskip

Notice that $L_{n+1}$ is isomorphic $\omega_{g,n}(\sum_{i=1}^n D_{i,n+1})$ as line bundles over $\oM_{g,n+1}=\overline{\mathcal{C}}_{g,n}$, respectively. This comes from the definition $L_{n+1}= \sigma_{n+1}^*(\omega_{g,n+1})$ and its comparison with $\omega_{g,n}$: they differ precisely at the nodal points of the fibers of $\pi$ (but this happens in codimension two in $\oM_{g,n+1}$) and at $D_{i,n+1}$ for each $i=1,\ldots,n$. This is useful, for instance, to compute $\kappa_0 = 2g-2+n$ as the Euler characteristic of the marked fibers of $\pi$ multiplied by $-1$. We have $\kappa_a = \pi^* \kappa_a +\psi_{n+1}^a$ and
$$\pi_*(\psi_1^{k_1} \ldots\psi_n^{k_n} \psi_{n+1}^{a+1}) = \psi_1^{k_1} \ldots\psi_n^{k_n} \kappa_a. $$

\medskip

Finally, another natural vector bundle on $\oM_{g,n}$ is the bundle whose fiber at $[C] \in \oM_{g,n}$ is the $g$-dimensional vector space of holomorphic differentials on $C$. More precisely the Hodge bundle $\mathbb{E}$ on $\oM_{g,n}$ is the rank $g$ vector bundle corresponding to the sheaf $\pi_*\omega_{g,n}$, whose fiber at $[C]$ is given by $H^0(C,\Omega_C^1)$. We denote by $\Lambda(s):= \sum_{i=0}^g s^i \lambda_i$ its Chern polynomial, where $\lambda_{i}=c_i(\mathbb{E}) \in H^{2i}(\oM_{g,n})$. Mumford proved that $\Lambda(s)\Lambda(-s)=1$.

\medskip

\begin{example}
On $\oM_{1,1}$ there is a canonical isomorphism between the Hodge bundle $\mathbb{E}$ and the relative cotangent line bundle $L_1$ given by $\omega \mapsto \omega|_{x_1}$ for a holomorphic differential $\omega$ on the marked elliptic stable curve $C$ (a nonzero regular differential on a marked elliptic curve does not vanish at any smooth point). This means that $\psi_1=\lambda_1$ on $\oM_{1,1}$. Moreover, by the first part of Remark \ref{remark:psi-pullback}, we have, on $\oM_{1,n}$,
\begin{equation}\label{eq:psi-on-Mbar1n}
\psi_i=\lambda_1 + \sum_{\substack{I\sqcup J = [n] \\ i\in I}} \tikz[baseline=-1mm]{\coordinate (A) at (0,0);\coordinate (B) at (12mm,0);\draw (A)--(B);\leg{A}{-60};\leg{A}{-120};\leg{B}{-60};\leg{B}{-120};\gg{0}{A};\gg{1}{B};\lab{A}{-90}{6.9mm}{...};\lab{A}{-90}{10mm}{\underbrace{\phantom{aaa}}_{\text{$I$}}};\lab{B}{-90}{6.9mm}{...};\lab{B}{-90}{10mm}{\underbrace{\phantom{aaa}}_{\text{$J$}}};}.
\end{equation}
We can use this to further compute $\kappa_1$ on $\oM_{1,1}$. By definition $\kappa_1= \pi_*\psi_2^2$ where $\pi$ forgets the second marked point on $\oM_{1,2}$. By the above formula
$$\psi_2 = \lambda_1 + D_{1,2}$$
and multiplying it by a further $\psi_2$ and remembering that $\psi_2 D_{1,2}=0$ we get $\psi_2^2 = \psi_2 \lambda_1$ on~$\oM_{1,2}$. This implies that $\kappa_1=\pi_*\psi_2^2=(\pi_* \psi_2)\lambda_1= \kappa_0 \lambda_1= \lambda_1$ on $\oM_{1,1}$.
\end{example}

\medskip

Recall from the theory of characteristic classes that $\Lambda(s)=\exp\left(\sum_{d\geq 1} (-1)^{d-1} (d-1)! s^d \ch_d(\mathbb{E})\right)$ where $\ch_d(\mathbb{E}) \in H^{2d}(\oM_{g,n})$ is the $d$-th Chern character of $\mathbb{E}$. Using the Grothendieck--Riemann-–Roch formula one can prove \cite{Zvo06} that 
\begin{equation}\label{eq:Mumford}
\ch_d(\mathbb{E}) = \frac{B_{d+1}}{(d+1)!} \left(\kappa_d - \sum_{i=1}^n \psi_i^d + \frac{1}{2} \sum_{\Gamma}\frac{1}{|\Aut(\Gamma)|} (\xi_\Gamma)_*\left(\frac{\psi_{e'}^d-(-\psi_{e''})^d}{\psi_{e'} + \psi_{e''}}\right)\right),\quad d\ge 1,
\end{equation}
where the sum runs over all stable graphs with one edge and labelling of the half-edge as $e=(e',e'')$ (the two possible labelling of the half-edges of $e$ are compensated by the factor of $\frac{1}{2}$ in front of the sum). In particular, we see that $\ch_d(\mbE)=0$ for even $d>0$, which also follows from the relation $\Lambda(s)\Lambda(-s)=1$. 

\medskip

\begin{example}
Let us apply equation \eqref{eq:Mumford} to compute $\lambda_1$ on $\oM_{1,1}$. Since $\oM_{1,1}$ has complex dimension $1$, we have $\Lambda(s)=1+s\lambda_1=\exp(s\cdot\ch_1(\mathbb{E}))= 1+ s\cdot \ch_1(\mathbb{E})$, which implies $\lambda_1=\ch_1(\mathbb{E})$. According to formula \eqref{eq:Mumford} we have $\ch_1(\mathbb{E})= \frac{B_2}{2}\left(\kappa_1 - \psi_1 + \frac{1}{2} \tikz[baseline=-1mm]{\coordinate (A) at (0,0);\lp{A}{0};\legm{A}{180}{1};\gg{0}{A};}\right)$. Now recall that $\kappa_1=\psi_1=\lambda_1$ on $\oM_{1,1}$ to obtain $\kappa_1=\psi_1=\lambda_1=\frac{1}{24} \tikz[baseline=-1mm]{\coordinate (A) at (0,0);\lp{A}{0};\legm{A}{180}{1};\gg{0}{A};}$ .
\end{example}

\medskip

On moduli spaces of stable curves of compact type, the forgetful morphism $\pi\colon\cM_{g,n+1}^\ct \to \cM_{g,n}^\ct$ is proper and psi, kappa and lambda classes can be defined either intrinsically or as restrictions of the corresponding classes from $\oM_{g,n}$. An important fact that we will use in what follows is that $\lambda_g$ vanishes when restricted to $\Delta_{\text{irr}} = \oM_{g,n} \setminus \cM_{g,n}^\ct$. Indeed, let $\xi\colon\oM_{g-1,n+2} \to \oM_{g,n}$ be the loop gluing map, then we can construct a nowhere vanishing morphism from $\xi^*\mathbb{E}$ to the trivial line bundle as follows: the fiber of the vector bundle $\xi^*\mathbb{E}$ at $[C]\in\oM_{g-1,n+2}$ is given by $\{\omega\in H^0(C,\Omega^1_C(x_{n+1}+x_{n+2}))|\res_{x_{n+1}}\omega+\res_{x_{n+2}}\omega=0\}$ and let us send $\omega$ from this space to $\res_{x_{n+2}}\omega$. The map $\xi$ is a degree $2$ cover of $\Delta_{\text{irr}}$, so this guarantees that the top Chern class of $\mathbb{E}|_{\Delta_{\text{irr}}}$ vanishes in cohomology with rational coefficients. This vanishing can be seen as a special case of a more general fact. For any stable graph $\Gamma$, the Hodge polynomial~$\Lambda(s)$ satisfies $\xi_\Gamma^* \Lambda(s) = \prod_{v\in V(\Gamma)} p_v^*\Lambda(s)$.

\medskip

\subsection{Strata algebra and tautological ring}
To study self-intersection of strata we are interested in their normal bundle. Notice first that the normal bundle in $\oM_{g,n}$ to the divisor $D$ corresponding to any stable graph with one edge is given by the tensor product of the tangent lines at each side of the node represented by that edge (moving away from the divisor $D$ means smoothing out the node, i.e. choosing a deformation $xy=t$ of the node $xy=0$ in the local model around the node). More in general, since the boundary divisor of $\oM_{g,n}$ has transverse intersections, for any stable graph $\Gamma$,
$$N_{\xi_\Gamma}= \bigoplus_{e \in E(\Gamma)} L^*_{e'} \otimes L^*_{e''},$$
where $e'$ and $e''$ are the half-edges forming the edge $e$ and $L_{e'}$ and $L_{e''}$ are the cotangent tautological line bundles at the marked points corresponding to the half-edges $e'$ and $e''$. This completes our control on intersection products of boundary strata $\oM_{^\Gamma}$. Indeed, in light of our discussion of the intersection of distinct irreducible components of the boundary divisor, together with the above normal bundle characterization and the fact that $c_1( L^*_{e'} \otimes L^*_{e''})=-(\psi_{e'}+\psi_{e''})$, we can write the intersection product of any two strata represented by dual graphs $\Gamma_1$ and $\Gamma_2$ in the following way. 

\medskip

For a stable graph $\Gamma$ and a subset $E_1\subset E(\Gamma)$, let us denote by $\Contr_{E_1}(\Gamma)$ the stable graph obtained by the contraction of all the edges from the set $E(\Gamma)\backslash E_1$. Consider two stable graphs $\Gamma_1,\Gamma_2\in G_{g,n}$. Let us denote by $S_{\Gamma_1,\Gamma_2}$ the set of isomorphism classes of tuples
$$
(\Gamma,E_1,E_2,\phi_1,\phi_2),\quad\text{where}\quad \begin{minipage}{10cm}
$\Gamma\in G_{g,n}$, $E_1,E_2\subset E(\Gamma)$, $E_1\cup E_2=E(\Gamma)$;\\
$\phi_i$, $i=1,2$, is an isomorphism $\Contr_{E_i}(\Gamma)\stackrel{\cong}{\to}\Gamma_i$.
\end{minipage}
$$
Note that the maps $\phi_1$ and $\phi_2$ are not uniquely determined by $\Gamma,E_1$, and $E_2$, because $\Gamma_1$ and $\Gamma_2$ can have nontrivial automorphisms. The map $\phi_i$ induces a bijection between the set $H(\Gamma_i)$ and a subset in $H(\Gamma)$. For $h\in H(\Gamma_i)$, we denote by $\phi_i^{-1}(h)$ the corresponding half-edge in $H(\Gamma)$. The map $\phi_i$ induces a surjective map $V(\Gamma)\to V(\Gamma_i)$. For $v\in V(\Gamma_i)$, we denote by $\phi_i^{-1}(v)$ the set of vertices in $V(\Gamma)$ mapping to $v$. We are now ready to write a desired product formula:
$$
(\xi_{\Gamma_1})_*(1) \cdot (\xi_{\Gamma_2})_*(1) = \sum_{(\Gamma,E_1,E_2,\phi_1,\phi_2)\in S_{\Gamma_1,\Gamma_2}}(\xi_{\Gamma})_*\left(\prod_{e \in E_1 \cap E_2}-(\psi_{e'}+\psi_{e''})\right).
$$
More in general, we can include the contribution of psi and kappa classes from the moduli at each vertex of the stable graphs to obtain the important formula
$$
(\xi_{\Gamma_1})_*(\gamma_1) \cdot (\xi_{\Gamma_2})_*(\gamma_2) = \sum_{(\Gamma,E_1,E_2,\phi_1,\phi_2)\in S_{\Gamma_1,\Gamma_2}} (\xi_{\Gamma})_*\left(\tgamma_1\tgamma_2\prod_{e \in E_1 \cap E_2}-(\psi_{e'}+\psi_{e''})\right),
$$
with 
$$
\gamma_i = \prod_{v\in V(\Gamma_i)} \prod_{a>0} p_{v}^*\kappa_{a}^{m_{a,v}}\prod_{h\in H(\Gamma_i)} \psi_h^{d_h} \in H^*(\oM_{\Gamma_i}),
$$ 
for $i=1,2$ and any $m_{a,v},d_h \in \mbZ_{\geq 0}$, and
$$
\tgamma_i=\prod_{v\in V(\Gamma_i)}\prod_{a>0}\lb\sum_{\tv\in\phi_i^{-1}(v)}p_{\tv}^*\kappa_a\rb^{m_{a,v}}\prod_{h\in H(\Gamma_i)}\psi_{\phi_i^{-1}(h)}^{d_h}\in H^*(\oM_{\Gamma}).
$$

\medskip

What cohomology classes can we reach with expressions like the ones appearing in the above formula, i.e. what is the $\mbQ$-vector subspace generated by classes of the form $(\xi_\Gamma)_* \gamma$, where $\Gamma$ is any stable graph and $\gamma = \prod_{v\in V(\Gamma)} \prod_{a>0} p_{v}^*\kappa_{a}^{m_{a,v}}\prod_{h\in H(\Gamma)} p_{v_h}^* \psi_h^{d_h} \in H^*(\oM_\Gamma)$, with $m_{a,v},d_h \in \mbZ_{\geq 0}$ such that $\sum_{a> 0} a m_{a,v} + \sum_{h \in H(v)} d_h \leq 3g(v)-3+n(v)$? It turns out that this system of vector subspaces, for any $g$ and $n$ in the stable range, can be characterized as the smallest system of subalgebras (with unit!) $R^*(\oM_{g,n}) \subset H^*(\oM_{g,n})$ closed under pushforward with respect to all gluing maps $\xi_\Gamma$ and forgetful maps $\pi$ forgetting any of the marked points. The subring~$R^*(\oM_{g,n})$ is called the tautological ring and its elements are called tautological classes. 

\medskip

The system of tautological rings is closed under pullback with respect to gluing and forgetful maps too and contains virtually any interesting geometrically natural class (to give a concrete example, the Hodge classes are tautological, but it is in fact highly nontrivial to exhibit or construct explicitly even a single non-tautological class). In genus $0$ it coincides with the full cohomology, but in higher genus it quickly becomes much smaller than cohomology. Notice that it cannot contain odd degree cohomology classes. By definition it is spanned by classes of the form described above, but the full system of relations among these generators is unknown (although several precise conjectures have been formulated). Any such relation is called a tautological relation.

\medskip

\begin{example}
In genus $0$ all tautological relations are generated by pullback or pushforward by the aforementioned WDVV relations in $H^0(\oM_{0,4})$. For instance, pulling back by the map forgetting $n$-further markings we obtain
\begin{equation}\label{eq:WDVVn+4}
\sum_{\substack {I\sqcup J =[n]}} \tikz[baseline=-1mm]{\coordinate (A) at (0,0);\coordinate (B) at (20mm,0);\draw (B)--(A);\legm{A}{180}{a};\legm{A}{90}{b};\leg{A}{-60};\leg{A}{-120};\gg{0}{A};\lab{A}{-90}{6.9mm}{...};\lab{A}{-90}{9mm}{I};\legm{B}{90}{c};\leg{B}{-60};\leg{B}{-120};\gg{0}{B};\lab{B}{-90}{6.9mm}{...};\lab{B}{-90}{9mm}{J};\legm{B}{0}{d};}= (b \leftrightarrow c),
\end{equation}
with $a,b,c,d$ any permutation of $n+1,n+2,n+3,n+4$.
\end{example}

\medskip


\section{Integrable systems of evolutionary PDEs}\label{section:PDEs}

In this lecture we will introduce commuting systems (hierarchies) of evolutionary PDEs, with a special focus on integrable systems of conservation laws. Our approach is mostly algebraic, as opposed to analytic and our focus is on symmetries and compatibility/commutativity of these systems, as opposed to constructing solutions. In view of the correspondence with intersection theory on moduli spaces of curves is, we propose a graphical representation for these compatibility equations, which will help us track down the corresponding tautological relations expressed in terms of stable graphs and strata algebra.

\medskip

\subsection{Evolutionary PDEs}

Evolutionary PDEs in one space and one time variables are systems of partial differential equations for $N$ unknown functions $u^\alpha=u^\alpha(x,t)$, $1\leq \alpha\leq N$, of two independent variables $x,t$ of the form
$$u^\alpha_t = Q^\alpha(u^*,u^*_x,u^*_{xx},\ldots),\qquad 1\leq\alpha\leq N,$$
where $u^\alpha_t:= \frac{\d u^\alpha}{\d t}$, $u^\alpha_x := \frac{\d u^\alpha}{\d x}$, $u^\alpha_{xx} := \frac{\d^2 u^\alpha}{(\d x)^2}$ and more in general $u^\alpha_k := \frac{\d^k u^\alpha}{(\d x)^k}$, and the symbol $*$ as a sub- or superscript means any possible value of the corresponding index, for instance $u^*$ means $u^\beta$ for $1\leq \beta\leq N$, while $u^\alpha_*$ means $u^\alpha_k$ for $k\geq 0$.

\medskip

Typically the right hand side $Q^\alpha= Q^\alpha(u^*,u^*_x,u^*_{xx},\ldots)$ is a differential polynomial, i.e. a polynomial in the variables $u^\alpha$ and their $x$-derivatives, but sometimes it's useful to include formal power series. In that case, however, we will sort the differential degree of each term of~$Q^\alpha$ by powers of a formal variable $\eps$, namely
$$
u^\alpha_t = Q^\alpha(u^*_*;\eps),\qquad 1\leq\alpha\leq N, \qquad \text{where} \qquad Q^\alpha \in \mbC[[u^*_*,\eps]], \qquad \deg Q^\alpha= d,
$$
where $d\in \mbZ_{\geq 0}$, with respect to the grading $\deg u^*_k = k$ and $\deg \eps=-1$. The reason behind it is that we want to be able to treat this equation order by order in $\eps$ so that, when we truncate at finite power of $\eps$ the right hand side contains only a finite number of $x$-derivatives, as a PDE should.

\medskip

\begin{example}
For $N=1$ we can suppress the index of the dependent variable $u^1=u$. The diffusion (heat) equation $u_t=u_{xx}$ is an example of evolutionary PDE with $d=2$. But if we want to study the Burgers equation $u_t= u u_x + u_{xx}$ instead, it might be convenient to realize it as a $d=1$ equation $u_t = u u_x + \eps u_{xx}$ in order to be able to deal with the nonlinear inviscid case first, for $\eps=0$.
\end{example}

\medskip

In fact the case $d=1$ is particularly suggestive because it can be seen as emerging from a ``long wave'' limit process consisting in applying to the $\eps$-independent case the equal rate space-time zoom $x\mapsto x/\eps$ and $t\mapsto t/\eps$. The $\eps=0$ limit is then a system of dispersionless evolutionary PDEs and the full series can be seen as a small dispersion expansion. Dispersion here refers to (the nonlinear analogue of) the fact that Fourier components of $u^\alpha(x)$ of different frequency propagate at different velocities. 

\medskip

\begin{example}
Our formal power series approach allows to accomodate difference-differential systems too. Indeed, the operator $\Lambda_\eps$ shifting a function $f(x)$ to a function $f(x+\eps)$ can be Taylor expanded as $\Lambda_\eps = \exp(\eps \d_x) = \sum_{k\geq 0}\frac{\eps^k \d_x^k}{k!}$. For instance the Toda lattice system
$$
\left\{\begin{array}{l} u^1_t=\frac{1}{\eps}\left(e^{u^2(x+\eps)}- e^{u^2(x)}\right),\\ 
u^2_t= \frac{1}{\eps}(u^1(x)-u^1(x-\eps)),\end{array}
\right. 
$$
belongs to the class of systems we are considering, with $d=1$, but can be seen as the evolution equations for a system of particles on a 1D lattice with spacing $\eps$ and first neighbours interaction.
\end{example}

\medskip

We can look at a system of evolutionary PDEs as a vector field on some infinite dimensional space in the same way as we interpret a system of ODEs $\frac{d y^\alpha}{dt} = Q^\alpha(y^*)$ as the local coordinate expression of the problem of finding the flow $t\mapsto y^\alpha(t)$ of the vector field $Q = Q^\alpha(y^*) \frac{\d}{\d y^\alpha}$ on a finite dimensional manifold with local coordinates $y^1,\ldots,y^N$. In this analogy the coordinate index $\alpha$ is replaced by the bi-index $(\alpha,x)$ where $\alpha$ is a discrete and $x$ a continuous index. Indeed the infinte dimensional manifold we are dealing with has coordinates $u^\alpha(x)$ and the flow of the vector field $Q$ gives rise to a curve whose coordinate description is $t\mapsto u^*(x,t)$. Of course, since in our approach $Q^\alpha$ is just a formal power series, the corresponding infinite dimensional manifold cannot be geometrically constructed in any meaningful way, but it's still heuristically useful to imagine there is one.

\medskip

\subsection{Compatibility and commutativity}

Given multiple systems of evolutionary PDEs, parametrized by some index set $\mathcal{I}$,
$$u^\alpha_{t_i} = Q^\alpha_i(u^*_*;\eps),\qquad 1\leq\alpha\leq N, \quad i\in \mathcal{I},$$
we can ask about their compatibility, i.e. existence of a common solution $u^*= u^*(x,t_*)$ depending simultaneously on all time variables $t_i$, for $i\in \mathcal{I}$. Compatibility is of course equivalent to the system of equations
$$\frac{\d^2 u^\alpha}{\d t_j \d t_i}= \frac{\d}{\d t_j} Q^\alpha_i(u^*_*;\eps) = \frac{\d}{\d t_i} Q^\alpha_j(u^*_*;\eps)=\frac{\d^2 u^\alpha}{\d t_j \d t_i},\qquad 1\leq \alpha\leq N,\quad i,j\in \mathcal{I},$$
which we can interpret as the vanishing of the vector fields Lie bracket $[Q_i,Q_j]$.

\medskip

\begin{example}
The Korteweg--de Vries equation is the most well known example of integrable nonlinear evolutionary PDE. It fits our framework with $N=1$ and $d=1$, namely
\begin{gather*}
u_t = u u_x + \frac{\eps^2}{12} u_{xxx}.
\end{gather*}
It turns out that the KdV equation is one member of an infinite hierarchy of commuting evolutionary PDEs, indexed by $\mathcal{I}=\mbZ_{\geq 0}$. We list the first few flows here, and the reader is invited to check commutativity:
\begin{equation*}
\begin{split}
u_{t_0}=&  u_x,\\
u_{t_1}=&  u u_x + \frac{\eps^2}{12} u_{xxx},\\
u_{t_2}=& \frac{1}{2}u^2u_x +\frac{\eps^2}{12}(2u_x u_{xx}+u u_{xxx})+\frac{\eps^4}{240}u_{xxxxx}.
\end{split}
\end{equation*}
\end{example}

\medskip

In fact the formalism can be made much more similar to the finite dimensional case if we assume that $u^\alpha(x)$ is a periodic function and we pass to its Fourier series $u^\alpha(x,t) = \sum_{a\in \mbZ} p^\alpha(t) e^{\mathrm{i}ax}$. This way we can use the Fourier coefficients $p^\alpha_a$, $1\leq \alpha\leq N$, $a\in\mbZ$ as an infinite but this time discrete, coordinate set. The evolution equations can then be written in the form
$$
(p^\alpha_a)_{t_i} = Q^\alpha_{i,a}(p^*_*,\eps), \qquad 1\leq \alpha\leq N,\quad i\in \mathcal{I}, \quad a\in\mbZ, \qquad \text{where} \qquad Q^\alpha_{i,a} \in \mbC[[p^*_*,\eps]].
$$
In this formalism it is much easier to compute the compatibility conditions between~$Q_i$ and~$Q_j$, and we obtain the familiar formula for the Lie bracket of vector fields,
\begin{equation}\label{eq:comm}
\sum_{b\in \mbZ}\left(\frac{\d Q^\alpha_{i,a}}{\d p^\mu_b} Q^\mu_{j,b} - \frac{\d Q^\alpha_{j,a}}{\d p^\mu_b} Q^\mu_{i,b} \right)=0, \qquad 1\leq \alpha\leq N,\quad i,j\in \mathcal{I}.
\end{equation}

\medskip

\begin{remark}
Naturally, the components $Q^\alpha_{i,a}\in \mbC[[p^*_*,\eps]]$ of the vector fields $Q_i$ cannot be arbitrary if we want to remember that $Q^\alpha_{i,a}$ is the $a$-th Fourier coefficient of $Q^\alpha_i\in \mbC[[u^*_*,\eps]]$, namely,
$$Q^\alpha_{i,a} = \sum_{\substack{a_1,\ldots,a_n\in\mbZ\\  g,n\geq 0}} Q^\alpha_{i,a,g;\alpha_1,\ldots,\alpha_n}(a_1,\ldots,a_n)\, \frac{p^{\alpha_1}_{a_1}\ldots p^{\alpha_n}_{a_n}}{n!} \eps^g \quad $$
with 
$$
Q^\alpha_{i,a,g;\alpha_1,\ldots,\alpha_n}(a_1,\ldots,a_n) = 
\left\{\begin{array}{ll} Q^\alpha_{i,g;\alpha_1,\ldots,\alpha_n}(a_1,\ldots,a_n),& \text{if $\sum a_i = a$}, \\ 0, & \text{otherwise},\end{array}\right.
$$
where $Q^\alpha_{i,g;\alpha_1,\ldots,\alpha_n}(a_1,\ldots,a_n)$ is a polynomial of degree $g+d$ in the $a_i$'s and independent of $a$.
\end{remark}

\medskip

We end this section with a graphical representation of the commutativity equation above. If we represent the $\sum a_i = a$ slice of the polynomial coefficient $Q^\alpha_{i,g;\alpha_1,\ldots,\alpha_n}(a_1,\ldots,a_n)$ by an $(n+2)$-valent vertex decorated as follows,
$$
Q^\alpha_{i,a,g;\alpha_1,\ldots,\alpha_n}(a_1,\ldots,a_n) \qquad \mapsto\qquad 
\tikz[baseline=-1mm]{\coordinate (A) at (0,0);\draw (A);\legmout{A}{180}{(\alpha,a)};\legm{A}{90}{i};\legin{A}{-60};\legin{A}{-120};\gg{g}{A};\lab{A}{-90}{6.9mm}{...};\lab{A}{-90}{9mm}{(\alpha_1,a_1),\ldots,(\alpha_n,a_n)};},
$$
and we represent summation over repeated indices by edges with matching half-edge decorations, then equation \eqref{eq:comm} can be represented as the system of equalities, valid for all $g,n\geq 0$, $1\leq \alpha,\alpha_1,\ldots,\alpha_n\leq N$, and $a,a_1,\ldots,a_n\in\mbZ$ with $\sum a_i = a$,
$$\sum_{\substack {g_1,g_2\geq 0\\ g_1+g_2=g\\ I\sqcup J =[n]}}
\sum_{\substack{1\leq \mu \leq N\\ b\in \mbZ}} \tikz[baseline=-1mm]{\coordinate (A) at (0,0);\coordinate (B) at (20mm,0);\draw[midarrow] (B)--(A);\legmout{A}{180}{(\alpha,a)};\legm{A}{90}{i};\legin{A}{-60};\legin{A}{-120};\gg{g_1}{A};\lab{A}{-90}{6.9mm}{...};\lab{A}{-90}{9mm}{(\alpha_I,a_I)};\legm{B}{90}{j};\legin{B}{-60};\legin{B}{-120};\gg{g_2}{B};\lab{B}{-90}{6.9mm}{...};\lab{B}{-90}{9mm}{(\alpha_J,a_J)};\lab{B}{165}{10.0mm}{(\mu,b)}} = (i \leftrightarrow j),$$
where the symbol $(i \leftrightarrow j)$ stands for the same expression as the left-hand side, but with~$i$ and~$j$ exchanged.

\medskip

\subsection{Systems of conservation laws}
In this section we focus on the special case of systems of the form
$$u^\alpha_{t_i} = \d_x P^\alpha_i(u^*_*;\eps),\qquad 1\leq\alpha\leq N, \quad i\in \mathcal{I},\qquad \text{where} \qquad P^\alpha \in \mbC[[u^*_*,\eps]], \qquad \deg P^\alpha= 0.$$
These are called systems of conservation laws because $\frac{d}{dt_i}\int u^\alpha dx = \int \d_x P^\alpha_i(u^*_*;\eps) = 0$ if there is no boundary contribution to the integral (e.g. if $u^\alpha$ is a function on $\mbR$ decaying fast enough together with its derivatives or if it is a function on $S^1$), so that $u^\alpha$ is the density of a conserved quantity.

\medskip

If we introduce Fourier coefficients $P^\alpha_{i,a}$ for the differential polynomials $P^\alpha_i$ in the same way as we did for $Q^\alpha_i$ in the previous section, we see that $Q^\alpha_{i,a}=\mathrm{i} a P^\alpha_{i,a}$ and the commutativity equations \eqref{eq:comm} become
\begin{equation}\label{eq:clcomm}
\sum_{b\in \mbZ}b\left(\frac{\d P^\alpha_{i,a}}{\d p^\mu_b} P^\mu_{j,b} - \frac{\d P^\alpha_{j,a}}{\d p^\mu_b} P^\mu_{i,b} \right)=0, \qquad 1\leq \alpha\leq N,\quad i,j\in \mathcal{I}.
\end{equation}
and, similarly
$$P^\alpha_{i,a} = \sum_{\substack{a_1,\ldots,a_n\in\mbZ\\ \sum a_i = a\\ g,n\geq 0}} P^\alpha_{i,g;\alpha_1,\ldots,\alpha_n}(a_1,\ldots,a_n)\, \frac{p^{\alpha_1}_{a_1}\ldots p^{\alpha_n}_{a_n}}{n!} \eps^g \quad $$
with
$$P^\alpha_{i,a,g;\alpha_1,\ldots,\alpha_n}(a_1,\ldots,a_n) = \left\{\begin{array}{ll} P^\alpha_{i,g;\alpha_1,\ldots,\alpha_n}(a_1,\ldots,a_n)& \text{ if } \sum a_i = a \\ 0 & \text{ otherwise }\end{array}\right.$$
and 
$P^\alpha_{i,g;\alpha_1,\ldots,\alpha_n}(a_1,\ldots,a_n)$ a polynomial of degree $g+d-1$ in the $a_i$'s.

\medskip

If we use the same graphical notation as in the last section, but this time applied to the coefficients $P^\alpha_{i,a,g;\alpha_1,\ldots,\alpha_n}(a_1,\ldots,a_n)$, which vanish unless $\sum a_i=a$, we get
\begin{equation}\label{eq:graphcomm}
\sum_{\substack {g_1,g_2\geq 0\\ g_1+g_2=g\\ I\sqcup J =[n]}} 
\sum_{\substack{1\leq \mu \leq N\\ b\in \mbZ}} 
b\left( \tikz[baseline=-1mm]{\coordinate (A) at (0,0);\coordinate (B) at (20mm,0);\draw[midarrow] (B)--(A);\legmout{A}{180}{(\alpha,a)};\legm{A}{90}{i};\legin{A}{-60};\legin{A}{-120};\gg{g_1}{A};\lab{A}{-90}{6.9mm}{...};\lab{A}{-90}{9mm}{(\alpha_I,a_I)};\legm{B}{90}{j};\legin{B}{-60};\legin{B}{-120};\gg{g_2}{B};\lab{B}{-90}{6.9mm}{...};\lab{B}{-90}{9mm}{(\alpha_J,a_J)};\lab{B}{165}{10.0mm}{(\mu,b)}}\right) = (i \leftrightarrow j).
\end{equation}

\medskip

\subsection{Dispersionless conservation laws}
Before trying to understand how, in general, the above equation is our bridge to intersection theory on the moduli space of curves, as a warm-up, let us consider the simpler case of dispersionless evolutionary PDEs. If our evolutionary PDEs have differential degree $d=1$ and we take the dispersionless limit $\eps=0$, we obtain a system of the form $u^\alpha_{t_i} = \d_x P^\alpha_i(u^*)$, where $P^\alpha_i(u^*)$ does not depend on $x$-derivatives of the $u^\alpha$ anymore. This means that the coefficients $P^\alpha_{i,0;\alpha_1,\ldots,\alpha_n}(a_1,\ldots,a_n)$ are constant in the $a_i$'s,
$$
P^\alpha_{i,a,0;\alpha_1,\ldots,\alpha_n}(a_1,\ldots,a_n) = \left\{\begin{array}{ll} P^\alpha_{i,0;\alpha_1,\ldots,\alpha_n},& \text{if $\sum a_i = a$}, \\ 0, & \text{otherwise},\end{array}\right.
$$
and
$$P^\alpha_i(u^*_*) = \sum_{n\geq 0}P^\alpha_{i,0;\alpha_1,\ldots,\alpha_n} \frac{u^{\alpha_1} \ldots u^{\alpha_n}}{n!} .$$
Now notice that in order to deduce an equation for the $P^\alpha_{i,0;\alpha_1,\ldots,\alpha_n}$ it is enough to consider the case only one of the $a_i$'s does not vanish. Let's then consider $n+1$ total ingoing legs with $a_1=\ldots=a_n=0$ and mark the last ingoing leg with $(\alpha_{n+1}=\beta,a_{n+1}=1)$. Then each summand in \eqref{eq:graphcomm} vanishes unless $n+1\in J$ and $b=1$ and we obtain the simpler graphical equation
\begin{equation}\label{eq:dlessgraphcomm}
\sum_{\substack {I\sqcup J =[n]}} 
\sum_{\substack{1\leq \mu \leq N}} 
\tikz[baseline=-1mm]{\coordinate (A) at (0,0);\coordinate (B) at (20mm,0);\draw[midarrow] (B)--(A);\legmout{A}{180}{\alpha};\legm{A}{90}{i};\legin{A}{-60};\legin{A}{-120};\gg{0}{A};\lab{A}{-90}{6.9mm}{...};\lab{A}{-90}{9mm}{\alpha_I};\legm{B}{90}{j};\legin{B}{-60};\legin{B}{-120};\gg{0}{B};\lab{B}{-90}{6.9mm}{...};\lab{B}{-90}{9mm}{\alpha_J};\legmin{B}{0}{\beta}:\lab{B}{165}{10.0mm}{\mu}}= (i \leftrightarrow j),
\end{equation}
where, now,
$$P^\alpha_{i,0;\alpha_1,\ldots,\alpha_n} \qquad \mapsto\qquad 
\tikz[baseline=-1mm]{\coordinate (A) at (0,0);\draw (A);\legmout{A}{180}{\alpha};\legm{A}{90}{i};\legin{A}{-60};\legin{A}{-120};\gg{0}{A};\lab{A}{-90}{6.9mm}{...};\lab{A}{-90}{9mm}{\alpha_1,\ldots,\alpha_n};}.$$
\begin{remark}
We could have started our discussion on commutativity and graphical representation of equations for systems of conservation laws with the disperionless case $u^\alpha_{t_i} = \d_x P^\alpha_i(u^*)$ by letting $$P^\alpha_i(u^*_*) =: \sum_{n\geq 0}P^\alpha_{i,0;\alpha_1,\ldots,\alpha_n} \frac{u^{\alpha_1} \ldots u^{\alpha_n}}{n!}$$ and easily deducing from $\frac{\d P^\alpha_i}{dt_j}=\frac{\d P^\alpha_j}{dt_i}$ the equations
$$\frac{\d P^\alpha_i}{\d u^\mu}\frac{\d P^\mu_j}{\d u^\gamma} =\frac{\d P^\alpha_j}{\d u^\mu}\frac{\d P^\mu_i}{\d u^\gamma}.$$
These are known as oriented WDVV equations, and their graphical representation is exactly~\eqref{eq:dlessgraphcomm}. The Fourier approach, however, is useful when considering the dispersive case.
\end{remark}

\medskip

\section{Moduli spaces of curves and the DR hierarchy}

Following our purely algebraic approach to systems of conservation laws, we will explain how to draw a very suggestive correspondence between the expressions of these PDEs and generating functions of intersection numbers on moduli space of curves so that properties of the systems of equations can be mapped to (tautological) relations in the cohomology of the moduli spaces. Even more interestingly, we will introduce a vast family of such correspondences, namely cohomological field theories and their generalizations, each of which produces a specific integrable hierarchy of PDEs possessing the universal properties controlled by the aforementioned tautological relations. This phenomenon can be distilled into the following one-liner: \emph{integrability is woven into the structure of intersection theory of $\oM_{g,n}$, cohomological field theories materialize it into specific integrable systems}.

\medskip

The final objective of this lecture is reviewing the construction and some of the properties of the double ramification (DR) hierarchy introduced by the second named author in \cite{Bur15-DR-hierarchy}, although we maintain a higher level of generality, in the spirit of \cite{ABLR21}. Useful references are the reviews \cite{Bur24,Ros17} and the papers \cite{BR16a, BDGR18,BDGR20,BGR19}. Along the way the notion of F-cohomological field theory will emerge as natural from the classical integrability viewpoint.

\subsection{Dispersionless conservations laws and WDVV relations}
Looking at equation \eqref{eq:dlessgraphcomm} we notice that we have encountered a very similar pictorial formula before, as the WDVV relations \eqref{eq:WDVVn+4} in $H^*(\oM_{0,n+4})$, only it was much less decorated: the graphs carried a similar genus $0$ labeling at the vertices, but the legs where marked just by an integer in $[n+4]$ and edges where not marked at all.

\medskip

But even more foundationally, graphs like the ones appearing in the WDVV relations \eqref{eq:WDVVn+4} represent cycles in the moduli space of curves, while those appearing in \eqref{eq:dlessgraphcomm} are (products and sums) of numbers of type $P^\alpha_{i,0;\alpha_1,\ldots,\alpha_n}$. So what is a natural way to transform a cycle into a number? Integrating some cohomology class of the right dimension over the cycle! In fact, if the cohomology class we integrate on the cycle depends on some extra parameters and these parameters are somehow associated to half-edges or vertices of the stable graphs, we might even account for the richer labeling.

\medskip

The idea we are after is exploiting cohomological relations in the moduli spaces which we know how to prove to produce interesting examples of integrable systems. 

\medskip

Given the linear nature of $H^*(\oM_{0,n+1})$, one way to do this is defining, for any $n\geq 2$, linear maps $c_{0,n+1}\colon V^*\otimes V^{\otimes n} \to H^*(\oM_{0,n+1})$, for vector spaces $V$ with a basis $\{e_\alpha \, |\, 1\leq \alpha\leq N\}$. Coherently with the view that that these vector ``insertions'' are attached to marked points, we require the maps $c_{0,n+1}$ to be $S_{n}$\emph{-equivariant} with respect to permutations of the last $n$ tensor factors of $V$ and the labeling of the last $n$ marked points on $\oM_{0,n+1}$. The reader might be wondering why we chose to attach a covector at the first marking, and the answer is that we see at the edge of the graphs appearing in \eqref{eq:dlessgraphcomm} that the corresponding covector index of one vertex will be saturated with a vector index from the other vertex. This way we produce the extra decoration at the half-edges of our graphs: a basis vector insertion $e_\alpha$ is represented by an ingoing half-edge with decoration $\alpha$, a basis covector insertion $e^\alpha$ is represented by an outgoing half-edge decoration $\alpha$.

\medskip

Summarizing, our working hypothesis now is
$$
P^\alpha_{(\beta,d),0;\alpha_1,\ldots,\alpha_n}= \int_{\oM_{0,n+2}} c_{0,n+2}(e^\alpha\otimes e_\beta\otimes e_{\alpha_1}\otimes\ldots\otimes e_{\alpha_n} ) \psi_2^d.
$$
Regarding the second marking, notice how we have decided to attach to it a power of the $\psi_2$ class, recalling that psi classes automatically come attached to markings. This accounts for the different type of labeling at that point, which now comes from the index set $\mathcal{I}=[N]\times \mbZ_{\geq 0}$. Here we should remark that we didn't really need to use the same vector space $V$ for the insertion at this special, unoriented edge/marking. A posteriori one motivation for this is that known examples of integrable systems of evolutionary PDEs most often have their flows labeled by $\mathcal{I}=[N]\times \mathbb{Z}_{\geq 0}$, where $N$ is the number of dependent variables in the equations.

\medskip

We still need to ensure that equation \eqref{eq:dlessgraphcomm} is satisfied by the above ``vertex'' contributions. Let $\gl=\xi_\Gamma$ with $\Gamma$ any stable tree with $n+1$ legs, two vertices with markings $\{1\}\sqcup I$ and $J$, respectively, $I\sqcup J=\{2,\ldots,n+1\}$, and one separating edge, then we require the following \emph{splitting axiom}:
$$
\gl^*\left(c_{0,n+1}(e^\alpha\otimes \otimes_{i\in [n]} e_{\alpha_i})\right) = c_{0,|I|+2}(e^\alpha\otimes \otimes_{i\in I} e_{\alpha_{i-1}}\otimes e_\mu)\times c_{0,|J|+1}(e^\mu\otimes \otimes_{j\in J} e_{\alpha_{j-1}}).
$$
With this assumption we easily obtain \eqref{eq:dlessgraphcomm} by integrating the class $c_{0,n+4}(e^\alpha\otimes e_{\beta_1} \otimes e_{\beta_2} \otimes e_\gamma \otimes \otimes_{i\in [n]} e_{\alpha_i})\psi_2^{d_1} \psi_3^{d_2}$ against equation \eqref{eq:WDVVn+4} with $a=1,b=2,c=3,d=4$: indeed trading the pushforward by $\gl=\xi_\Gamma$ in the cycle with the pullback by $\gl$ of the above class and using the aforementioned splitting axiom we obtain the result.

\medskip

We see how the same, universal WDVV relation will guarantee commutativity for all hierarchies of dispersionless conservation law of the form
$$
u^\alpha_{t^\beta_d} = \d_x P^\alpha_{\beta,d}(u^*), \qquad 1\leq \alpha,\beta\leq N, \quad d\geq 0,
$$
with
$$
P^\alpha_{\beta,d} = \sum_{n\geq 1} \int_{\oM_{0,n+2}} c_{0,n+2}(e^\alpha\otimes e_\beta\otimes \otimes_{i=1}^n e_{\alpha_i}) \psi_2^d \frac{1}{n!}\prod_{i=1}^n u^{\alpha_i},
$$
where, for any $n\geq 2$, $c_{0,n+1}$ satisfies $S_n$-equivariance and the splitting axiom.

\medskip

\subsection{Dispersive conservation laws and the DR cycle}

Now we are ready for the general case of dispersive systems of conservation laws, where commutativity is encoded into equations~\eqref{eq:graphcomm}. We can expect that an all genus extension of the classes $c_{0,n+1}$ considered in the previous section will be needed and that it will take care of the vector half-edge labeling of the stable graphs therein. But before doing that, where do we find, in the cohomology of $\oM_{g,n+3}$, a relation corresponding to equation \eqref{eq:graphcomm}? As it turns out, this looks exactly like the pullback of the WDVV relations to the space of relative stable maps to (rubber) $\mbP^1$.

\medskip

The moduli space $\oM_{g,A} (\mbP^1;0,\infty)^\sim$ of stable maps to rubber $\mathbb{P}^1$ relative to $0$ and $\infty$, as a set, contains stable maps from nodal curves to $\mbP^1$ with ramif\/ication prof\/ile (orders of poles and zeros) given by $A=(a_1,\dots,a_n)\in \mathbb{Z}^n$, where $\sum\limits_{i=1}^n a_i= 0$. Here ``rubber'' means that we consider maps up to the $\mbC^*$-action in the target $\mbP^1$ and a~positive/negative coef\/f\/i\-cient~$a_i$ indicates a pole/zero at the $i$-th marked point, while $a_i=0$ just indicates an internal marked point (that is not a zero or pole). It's compactification and boundary stratification were studied in \cite{Li01,Li02}, where the author proved that $\oM_{g,A} (\mbP^1;0,\infty)^\sim$ can be endowed with the structure of a separated and proper Deligne--Mumford stack with a virtual fundamental class of homological degree $2(2g-3+n)$.

\medskip

A stable map from a smooth curve $C$ to rubber $\mbP^1$ is a branched cover with fixed branch profile over $0$ and $\infty$ dictated by $A$. Here is a picture for when $A\neq (0,\ldots,0)$.
\begin{figure}[H]
\begin{center}
\raisebox{-0.45\height}{\includegraphics[width=0.3\textwidth]{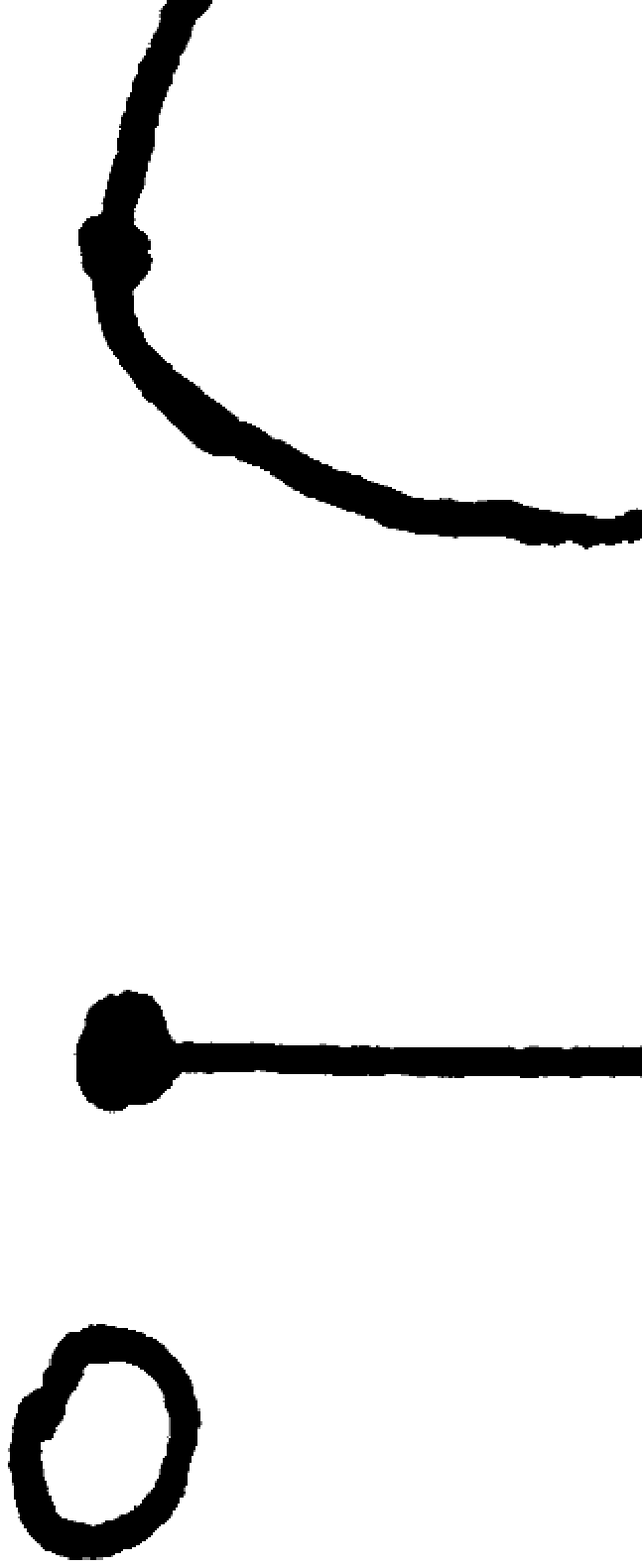}}
\end{center}
\end{figure}

When $A=(0,\ldots,0)$ the whole curve maps to a point in $\mbP^1$. At the boundary of $\oM_{g,A} (\mbP^1;0,\infty)^\sim$ we find stable maps from nodal curves $C$ to a chain of rubber $\mbP^1$'s attached to each other by a node (with $0$ and $\infty$ at opposite ends of the chain, so the target rubber $\mbP^1$ of the stable map degenerates together with the genus $g$ source). Over each irreducible rubber $\mbP^1$ in the chain there is a possibly disconnected stable branched cover with fixed ramification profile over the two nodes, or a node and either $0$ or $\infty$, from some of the irreducible components of the nodal genus~$g$ curve $C$. Over $0$ and $\infty$ the ramification profile is the one dictated by $A$, while the preimage of the nodes of the chain is a subset of the nodes of $C$ with opposite ramification multiplicity at each side of these nodes. The other nodes of $C$ are mapped to interior points in the chain of $\mbP^1$'s. Notice that there might be contracted irreducible components, mapping to a point of~$\mbP^1$, whose markings carry ramification multiplicity $0$. A stable map is one such nodal cover where at least one of the connected components over each $\mbP^1$ must be a stable component, i.e. different from the map $z^a$, for some positive integer $a$, from an otherwise unmarked $\mbP^1$ (which we refer to as unstable components). Notice that unstable components are not affected by the $\mbC^*$-action in the target with respect to which we take a quotient, while stable components are, and a covering map over a fixed $\mbP^1$ in the target chain with $r$ stable connected components has a residual $(\mbC^*)^{r-1}$-symmetry given by relative rescaling of the covering map. Here is a picture (we have omitted internal marked points or contracted components, but both are possible):\\
\begin{figure}[H]
\begin{center}
\raisebox{-0.45\height}{\includegraphics[width=0.5\textwidth]{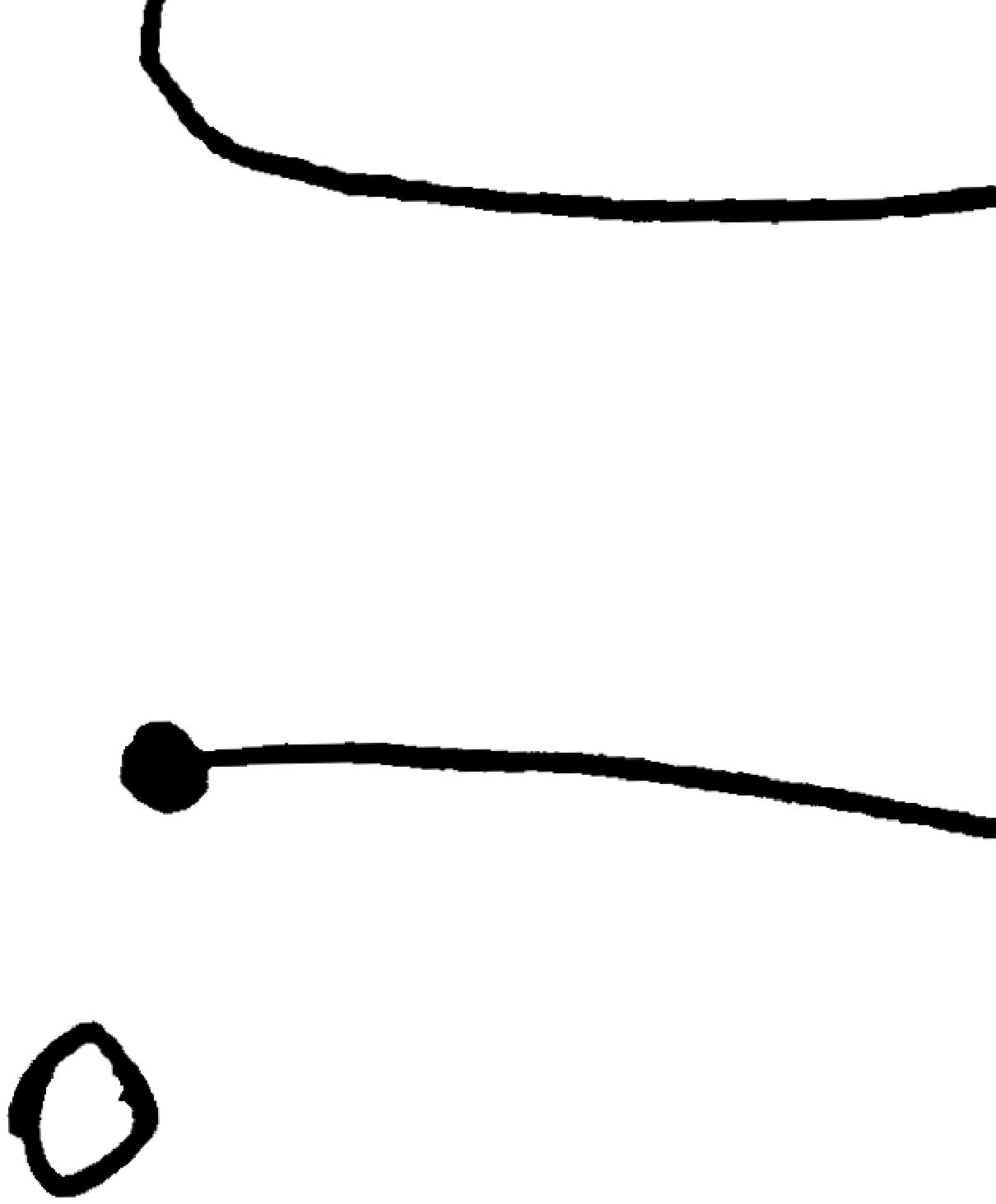}}
\end{center}
\end{figure}

\medskip

Although $\oM_{g,A} (\mbP^1;0,\infty)^\sim$ is not equidimensional (in particular because of the aforementioned contracted components), its virtual dimension is $2g-3+n$, and curves to a chain of $\mbP^1$'s of length $m$ appear in a boundary stratum of virtual codimension $m$.

\medskip

Let $p\colon\oM_{g,A} (\mbP^1;0,\infty)^\sim \to \oM_{g,n}$ be the map that only remembers the source curve, after stabilization. Then we define the double ramification cycle as
$$
\DR_g(A) := p_*([\oM_{g,A} (\mbP^1;0,\infty)^\sim ]^\vir).
$$
We view the DR cycle as a cohomology class in $H^{2g}(\oM_{g,n})$ dual to the homology class represented by some natural compactif\/ication of the locus, inside $\mathcal{M}_{g,n}$, formed by complex curves with marked points $x_1,\dots,x_n$ such that $\sum a_i x_i$ is the divisor of the zeros and poles of a~meromorphic function. We will adopt the convention $\DR_g(a_1,\ldots,a_n):=0$ if $\sum a_i\ne 0$.

\medskip

The class $\DR_g(a_1,\dots,a_n)$ is tautological and a (non-homogeneous) polynomial in the $a_i$'s formed by monomials of all even degrees from $0$ to $2g$. It's easy to see that it has the property $\pi^*\DR_g(a_1\ldots,a_n) = \DR_g(a_1\ldots,a_n,0)$, where $\pi$ forgets the last marking.

\medskip

\begin{remark} The above definition of the DR cycle is quite technical and it involves virtual classes, but it turns our it is a tautological class and we are lucky enough that people computed its expression in terms of generators of the tautological ring, so that we can work with that in actual computations. Here, we recall Pixton's very explicit construction of a nonhomogeneous cohomology class on $\oM_{g,n}$, with nontrivial terms in degree $0,2,4,\ldots,2g$. By a result of~\cite{JPPZ17}, the degree~$2g$ part of this class coincides with the double ramification cycle.

\medskip

Let $A=(a_1,\ldots,a_n)$ be a tuple of integers satisfying $\sum_{i=1}^n a_i=0$. Let $\Gamma$ be a stable graph and $r\ge 1$. A {\it weighting mod $r$} of $\Gamma$ is a function
$$
w\colon H(\Gamma)\to \{0,\ldots,r-1\}
$$
that satisfies the following three properties:
\begin{enumerate}
\item for any leg $i\in L(\Gamma)$, we have $w(i)=a_i\mod r$;
\item for any edge $e\in E(\Gamma)$, we have $w(e')+w(e'')=0\mod r$;
\item for any vertex $v\in V(\Gamma)$, we have $\sum_{h\in H(\Gamma),\,v_h=v}w(h)=0\mod r$.
\end{enumerate}
Denote by $W_{\Gamma,r}$ the set of weightings mod $r$ of $\Gamma$. We have $|W_{\Gamma,r}|=r^{h^1(\Gamma)}$.

\medskip

We denote by $P_g^{d,r}(A)\in H^{2d}(\oM_{g,n})$ the degree $2d$ component of the cohomology class 
\begin{equation}\label{eq:Pixton}
\sum_{\Gamma\in G_{g,n}}\sum_{w\in W_{\Gamma,r}}\frac{1}{|\Aut(\Gamma)|}\frac{1}{r^{h^1(\Gamma)}}\xi_{\Gamma*}\left[\prod_{i=1}^n\exp(a_i^2\psi_{i})\prod_{e\in E(\Gamma)}\frac{1-\exp\left(-w(e')w(e'')(\psi_{e'}+\psi_{e''})\right)}{\psi_{e'}+\psi_{e''}}\right]
\end{equation}
in $H^*(\oM_{g,n})$. Note that the factor $\frac{1-\exp\left(-w(e')w(e'')(\psi_{e'}+\psi_{e''})\right)}{\psi_{e'}+\psi_{e''}}$ is well defined since the denominator formally divides the numerator. In~\cite{JPPZ17}, the authors proved that for fixed $g,A$, and $d$ the class $P_g^{d,r}$ is polynomial in~$r$ for all sufficiently large~$r$. Denote by $P_g^d(A)$ the constant term of the associated polynomial in~$r$.

\medskip

In~\cite{JPPZ17}, the authors proved that 
$$
\DR_g(A)=2^{-g}P^g_g(A).
$$
In~\cite{CJ18}, the authors proved that the class $P_g^d(A)$ vanishes for $d>g$. Moreover it is easy to see from considerations on virtual fundamental classes of the space of degree $0$ stable maps that $\DR_g(0,\ldots,0)=(-1)^g\lambda_g$. 
\end{remark}

\medskip

\begin{remark}
Pixton's formula gives, in particular, a formula of Hain \cite{Hai11} for the restriction of the DR cycle to the moduli space of curves of compact type. Indeed, in order to compute the restriction of the class $P^j_g(A)$ to $\cM_{g,n}^\ct$, it is enough to disregard all terms in \eqref{eq:Pixton} whose underlying graph is not a tree, and we obtain
$$
\left.P^j_g(A)\right|_{\cM^\ct_{g,n}}=\left.\frac{2^j}{j!}\theta_g(A)^j\right|_{\cM^\ct_{g,n}}, \quad \text{with }\ \theta_g(A):= -\frac{1}{4} \sum_{J\subset \{1,\dots,n\}} \sum_{h=0}^g a_{J}^2 \delta^J_h,
$$
where
\renewcommand{\gg}[2]{\fill[color=white] (#2) circle(3.5mm) node {\color{black}$\substack{#1}$}; \draw (#2) circle (3.5mm)}
\begin{gather*}a_J := \sum_{j\in J} a_j , \qquad\delta^J_h: = \tikz[baseline=-1mm]{\coordinate (A) at (0,0);\coordinate (B) at (12mm,0);\draw (B)--(A);\leg{A}{-60};\leg{A}{-120};\gg{h}{A};\lab{A}{-90}{6.9mm}{...};\lab{A}{-90}{9mm}{J};\leg{B}{-60};\leg{B}{-120};\gg{g-h}{B};\lab{B}{-90}{6.9mm}{...};\lab{B}{-90}{9mm}{[n]\setminus J};}, \qquad \delta^{\{i\}}_0=\delta_g^{[n]\backslash\{i\}}: = - \psi_i.
\end{gather*}
\renewcommand{\gg}[2]{\fill[color=white] (#2) circle(2.5mm) node {\color{black}$\substack{#1}$}; \draw (#2) circle (2.5mm)}This gives Hain's formula
\begin{gather}\label{eq:Hain}
H^{2g}\big(\mathcal{M}_{g,n}^{\text{ct}}\big)\ni \DR_g(a_1,\dots,a_n)\big|_{\mathcal{M}_{g,n}^{\text{ct}}} = \frac{1}{g!}\left( -\frac{1}{4} \sum_{J\subset \{1,\dots,n\}} \sum_{h=0}^g a_{J}^2 \delta^J_h \right)^{g}.
\end{gather}
We see that $\DR_g(a_1,\dots,a_n)\big|_{\mathcal{M}_{g,n}^{\text{ct}}}$ is a homogeneous polynomial of degree $2g$ in the $a_i$'s. In particular, we have $\DR_g(0,\ldots,0)\big|_{\mathcal{M}_{g,n}^{\text{ct}}} =0$.
\end{remark}

\medskip 

Consider $\oM_{g,(0,0,A)} (\mbP^1;0,\infty)^\sim$, where the first two markings are internal marked points (label them $x$ and $y$, for instance). This space is endowed with a morphism~$\mathrm{br}$ to $\oM_{0,4}$ sending a stable map from a genus $g$ marked curve $C$ to its image $\mbP^1$ with the four markings given by $0, \infty,$ and the image of the first two (interior) markings of $C$. We can pullback the WDVV relation corresponding to exchanging the images of $x$ and $y$, from $H^2(\oM_{0,4})$ to $H^2(\oM_{g,(0,0,A)} (\mbP^1;0,\infty)^\sim)$ and then (virtual) pushforward it to $\oM_{g,n+2}$. All components of $\mathrm{br}^*(\text{WDVV})$ are loci whose generic elements are stable curves to a chain of exactly two $\mbP^1$'s and pushing forward to $\oM_{g,n+2}$ kills all contributions from covers with more than one stable component over each $\mbP^1$, because of the residual $\mbC^*$-symmetry we mentioned above. So nontrivial contributions come from the loci of maps of the following form (we omit drawing unstable components)

\medskip

\begin{figure}[H]
\begin{center}
\raisebox{-0.45\height}{\includegraphics[width=0.3\textwidth]{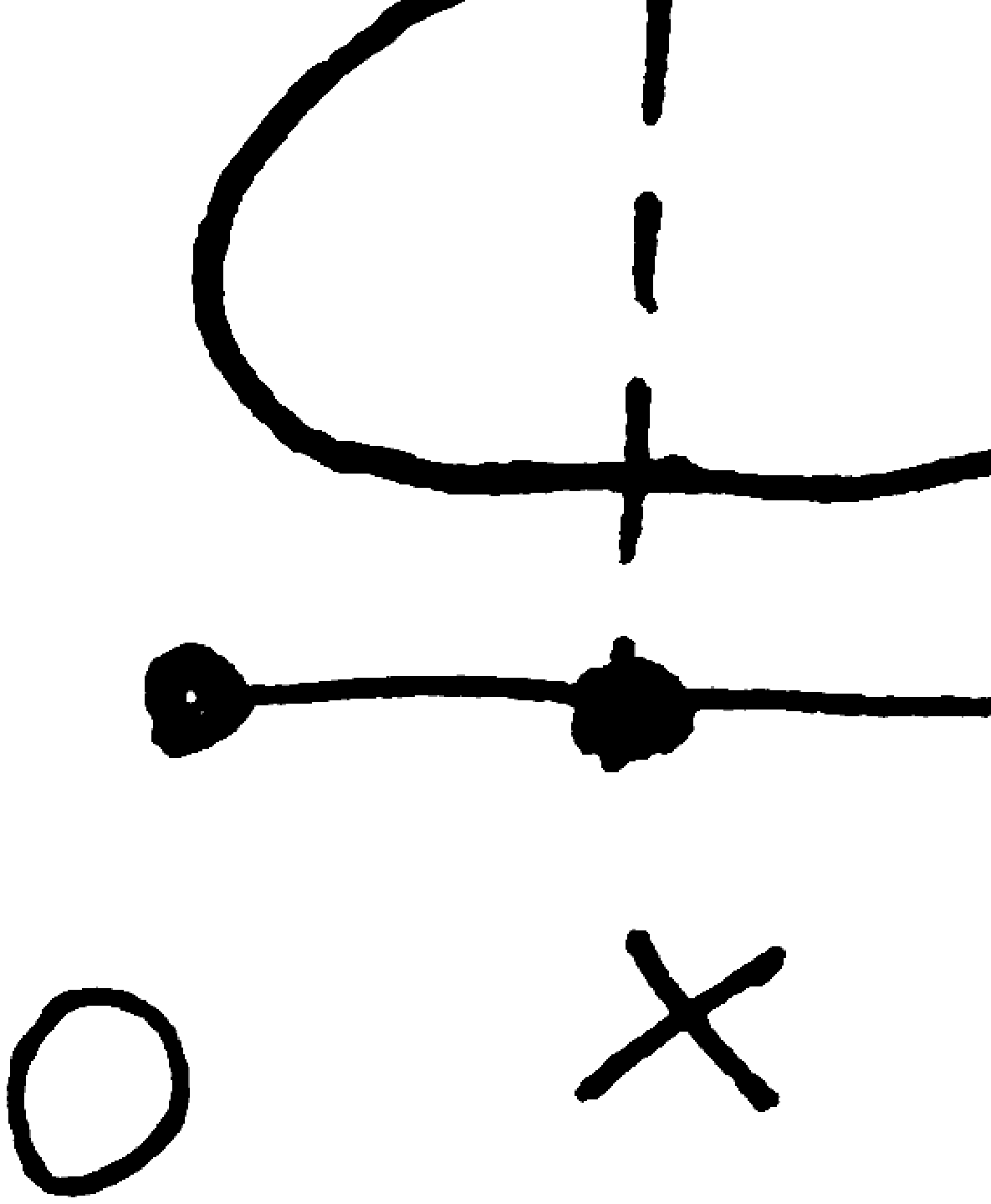}}
\end{center}
\end{figure}
\noindent together with the same maps where $x$ and $y$ are exchanged. After pushforward to $\oM_{g,n+2}$, what we obtain can be described more precisely as follows \cite{BSSZ15}. For $b_1,\ldots,b_p\in\mbZ_{>0}$ and $I\sqcup J=\{1,\ldots,n\}$, let us denote by $\DR_{g_1}(0_{1},A_I,-b_1, \dots, -b_p)\boxtimes \DR_{g_2}(0_{2},A_J, b_1, \dots, b_p)$ the cycle in~$\oM_{g_1+g_2+p-1,n+2}$ obtained by gluing the two double rami\-f\/i\-ca\-tion cycles at the marked points labeled by the positive integers~$b_1,\dots,b_p$. Here we denote by $A_I$ the tuple of integers formed by $a_i$, $i\in I$, and~$0_{1}$ and~$0_2$ indicate a coef\/f\/icient $0$ at the marked point~$1$ and~$2$. Then
\begin{gather}
\sum \frac{\prod\limits_{i=1}^p b_i}{p!}
\DR_{g_1}(0_{1},A_I,-b_1, \dots, -b_p)\boxtimes \DR_{g_2}(0_{2},A_J, b_1, \dots, b_p)\nonumber\\
\qquad{} - \sum \frac{\prod\limits_{i=1}^p b_i}{p!}
\DR_{g_1}(0_{2},A_I,-b_1, \dots, -b_p) \boxtimes \DR_{g_2}(0_{1},A_J, b_1, \dots, b_p)=0.\label{eq:DRcomm}
\end{gather}
The sums are over $I$, $J$, $p>0$ $b_1>0,\dots,b_p>0$, $g_1\ge 0$, $g_2\ge 0$, satisfying $g_1+g_2+p-1=g$. The factors $\frac{\prod\limits_{i=1}^p b_i}{p!}$ correspond to the degree of the map $\mathrm{br}$, which comes from the fact that there are $b$ ways of undegenerating a node with ramification number $b$.

\medskip

We have managed to upgrade the genus $0$ WDVV relations to relations in all genera and the idea, roughly speaking, was restricting to a cycle in $\oM_{g,n+2}$ endowed with a map to $\oM_{0,4}$ so that we could pull them back. If we wanted to find a cohomological counterpart to equations~\eqref{eq:graphcomm}, though, we are not there yet. If we interpret the vertices of the graphs appearing in \eqref{eq:graphcomm} as DR cycles, then we see they should be glued along one single edge, but equations \eqref{eq:DRcomm} have two DR cycles attached at any number $p$ of edges. But we already encountered the right tool to solve this problem: intersect this relation with the top Hodge class $\lambda_g$. Because $\lambda_g|_{\Delta_\mathrm{irr}}=0$, this has the effect of killing all terms with $p>1$ in \eqref{eq:DRcomm}, and, trading two sums over positive~$b$ for one sum over $b\in\mbZ$, we get
\begin{gather}
\sum_{\substack{b\in\mbZ \\ g_1,g_2\geq 0\\ g_1+g_2=g \\ I\sqcup J = \{1,\ldots,n\}} } b\ 
\lambda_{g_1}\DR_{g_1}(0_{1},A_I,-b)\boxtimes \lambda_{g_2}\DR_{g_2}(0_{2},A_J, b) =0.\label{eq:DRcommct}
\end{gather}
If we work with $n+3$ marked points $(0_1,0_2,a,a_1,\ldots,a_n)$ instead, singling out the third point, with ramification profile $a$, rewriting equation \eqref{eq:DRcommct} we get two terms, one where $a$ lies in the same component as $0_1$ and one where it lies in the same component as $0_2$. We then get exactly the cycle analogue of equation \eqref{eq:graphcomm}, with DR cycle classes at each vertex:
\renewcommand{\gg}[2]{\fill[color=white] (#2) circle(4.3mm) node {\color{black}$\substack{#1}$}; \draw (#2) circle (4.3mm)}
\begin{equation}\label{eq:graphDRcomm}
\sum_{\substack {g_1,g_2\geq 0\\ g_1+g_2=g\\ I\sqcup J =\{1,\ldots,n\}}} 
\sum_{\substack{ b\in \mbZ}} \ 
b\left( \tikz[baseline=-1mm]{\coordinate (A) at (0,0);\coordinate (B) at (20mm,0);\draw (B)--(A);\legm{A}{180}{a};\legm{A}{90}{0_1};\leg{A}{-60};\leg{A}{-120};\gg{\DR_{g_1}}{A};\lab{A}{-90}{6.9mm}{...};\lab{A}{-90}{9mm}{a_I};\legm{B}{90}{0_2};\leg{B}{-60};\leg{B}{-120};\gg{\DR_{g_2}}{B};\lab{B}{-90}{6.9mm}{...};\lab{B}{-90}{9mm}{a_J};\lab{B}{165}{10.0mm}{b}}\right) = (0_1 \leftrightarrow 0_2).
\end{equation}
\renewcommand{\gg}[2]{\fill[color=white] (#2) circle(2.5mm) node {\color{black}$\substack{#1}$}; \draw (#2) circle (2.5mm)}
Note that, because of the aforementioned vanishing properties of $\lambda_g$, we can actually work on~$\cM_{g,n}^{\mathrm{ct}}$ and Hain's formula~\eqref{eq:Hain} is enough for our purpose.

\medskip

\subsection{ct-F-CohFTs and integrable systems of conservation laws}

Inspired by what we saw in the $g=0$ case and with the aim of enriching equation \eqref{eq:graphDRcomm} with directed half-edges and Greek indices, we define a compact type F-cohomological field theory (ct-F-CohFT) as a family of linear maps $c_{g,n+1}\colon V^*\otimes V^{\otimes n} \to H^\even(\cM^{\ct}_{g,n+1})$ satisfying $S_n$-equivariance (in the same sense as we have seen in genus $0$) and the splitting axiom
$$
\gl^* c_{g,n+1}(e^\alpha\otimes \otimes_{i\in [n]} e_{\alpha_i}) = c_{g_1,|I|+2}(e^\alpha\otimes \otimes_{i\in I} e_{\alpha_{i-1}}\otimes e_\mu)\times c_{g_2,|J|+1}(e^\mu\otimes \otimes_{j\in J} e_{\alpha_{j-1}}),
$$
with respect to $\gl=\xi_\Gamma$ with $\Gamma$ any stable graph with $n+1$ legs, two vertices with genus~$g_1$ and~$g_2$ and markings $\{1\}\sqcup I$ and $J$, respectively, $g_1+g_2=g$ and $I\sqcup J=\{2,\ldots,n+1\}$, and one separating edge. F-CohFTs were first introduced in \cite{BR21}.

\medskip

Suppose that the classes of our ct-F-CohFT can be extended to $\oM_{g,n}$, then an arbitrary extension followed by the multiplication by $\lambda_g$ provides us with classes $\lambda_g c_{g,n+1}(e^\alpha\otimes \otimes_{j=1}^n e_{\alpha_j})$. The result is independent of the chosen extension, because if $\omega|_{\cM_{g,n}^\ct}=0$, then by the long exact sequence
$$
\cdots\to H^a(\oM_{g,n},\cM_{g,n}^\ct)\to H^a(\oM_{g,n})\to H^a(\cM^\ct_{g,n})\to H^{a+1}(\oM_{g,n},\cM_{g,n}^\ct)\to\cdots
$$
the class $\omega$ belongs to the image of the group $H^a(\oM_{g,n},\cM_{g,n}^\ct)$, which is isomorphic to $H_{2(3g-3+n)-a}(\Delta_{\text{irr}})$. Since $\lambda_g|_{\Delta_{\text{irr}}}=0$, we obtain $\lambda_g\omega=0$, as required.

\medskip

Then, upon defining
\begin{gather}\label{eq:DR hierarchy}
\frac{\d u^\alpha}{\d t^\beta_d}=\d_x P^\alpha_{\beta,d},\qquad 1\le \alpha,\beta\le N,\quad d\ge 0,
\end{gather}
where
\begin{equation}\label{eq:P-polynomials for DR hierarchy}
P^\alpha_{\beta,d}:=\sum_{\substack{g\ge 0,\,n\ge 1\\k_1,\ldots,k_n\geq 0\\\sum_{j=1}^n k_j=2g}} \frac{\eps^{2g}}{n!} \Coef_{(a_1)^{k_1}\ldots(a_n)^{k_n}} \left(\int_{\DR_g(-\sum_{j=1}^n a_j,0,a_1,\ldots,a_n)} \hspace{-2.3cm}\psi_2^d \lambda_g c_{g,n+2}(e^\alpha\otimes e_\beta\otimes \otimes_{j=1}^n e_{\alpha_j}) \right)\prod_{j=1}^n u^{\alpha_j}_{k_j},
\end{equation}
we obtain an infinite system of commuting evolutionary conservation laws, thanks to the relations \eqref{eq:graphDRcomm} for any ct-F-CohFT $c_{g,n+1}$. Indeed, commutativity of two of these flows is equivalent to the intersection of relations \eqref{eq:graphDRcomm} with the class $c_{g,n+3}(e^\alpha\otimes e_{\beta_1} \otimes e_{\beta_2} \otimes \otimes_{i=4}^{n+3}  e_{\alpha_i})$ (the reader should check it, using the axioms of ct-F-CohFT).

\medskip

\begin{remark}
If we don't assume that the classes of our ct-F-CohFT can be extended to $\oM_{g,n}$, then we propose the following modification of the above construction. Since $\lambda_g|_{\Delta_{\text{irr}}}=0$, by the exact sequence
$$
\cdots\to H^i_\cc(\cM_{g,n}^\ct)\stackrel{\iota}{\to} H^i_\cc(\oM_{g,n})\to H^i_\cc(\Delta_{\text{irr}})\to H^{i+1}_\cc(\cM_{g,n}^\ct)\to\cdots,
$$
where by $H^*_\cc(X)$ we denote the cohomology groups with compact supports, there exists a class $\tlambda_g\in H^{2g}_\cc(\cM_{g,n}^\ct)$ such that $\iota(\tlambda_g)=\lambda_g$. If $\omega\in H^*(\cM_{g,n}^\ct)$ is a class of our ct-F-CohFT, then the product~$\tlambda_g\omega$ belongs to the group $H^*_\cc(\cM_{g,n}^\ct)$ and we can define $\lambda_g\omega:=\iota(\tlambda_g\omega)$. We believe that appropriately chosen classes $\tlambda_g$ satisfy natural splitting properties such that the above construction of DR hierarchy produces a system of commuting flows.
\end{remark}

\medskip

\begin{example}\label{example:KdV}
The simplest example of ct-F-CohFT has $V=\mbC$ and $c_{g,n+1}(e^1\otimes e_1^{\otimes n})=1$ for all~$g,n$. Indeed the axioms of ct-F-CohFT are clearly trivially satisfied. The corresponding double ramification hierarchy was computed in \cite{Bur15-DR-hierarchy} using Hain's formula, and coincides with the KdV hierarchy. Let's try and reproduce that computation: we will compute the differential polynomials $P_d=P^1_{1,d}$. The relevant intersection numbers, according to formula~\eqref{eq:P-polynomials for DR hierarchy} above, are
$$
\int_{\DR_g(-\sum_{j=1}^n a_j,0,a_1,\ldots,a_n)} \lambda_g \psi_2^d.
$$
The cohomological degree of the integrand $\lambda_g \psi_2^d$ is $2(g+d)$ and we integrate it over the DR cycle, which has homological degree $2(2g-1+n)$, so we need $n=d-g+1$ to have a nontrivial contribution. When $g=0$ we get in particular $n=d+1$ and using $\DR_0(-\sum_{j=1}^n a_j,0,a_1,\ldots,a_n)=1$ and $\int_{\oM_{0,d+3}} \psi_2^d = 1$ we obtain the dispersionless limit of the hierarchy as
$$u_{t_d}=\frac{u^d}{d!} u_x + O(\eps^2), \qquad d\geq 0.$$
In general, when $d=0$, we must have $n=1$ and $g=0$. This gives
$$
u_{t_0} = u_x.
$$
When $d=1$, instead, there are two cases: we already considered $g=0$, $n=2$, but we also have $g=1$, $n=1$. We need to use Hain's formula to obtain
$$
\int_{\DR_1(-a,0,a)} \lambda_1 \psi_2 = 2 \int_{\DR_1(-a,a)} \lambda_1 = 2a^2 \int_{\oM_{1,2}}\delta_0^{\{1,2\}} \lambda_1 = 2 a^2\int_{\oM_{1,1}} \lambda_1=\frac{a^2}{12}.
$$
Here we first used that $\pi_* \psi_2 =\kappa_0= 2$ for $\pi\colon\oM_{1,3}\to \oM_{1,2}$ the map forgetting the second marking, then Hain's formula $\left.\DR_1(-a,a)\right|_{\cM_{1,2}^\ct}=\frac{a^2}{2}(\psi_1+\psi_2)$, which in turn is expressed through formula \eqref{eq:psi-on-Mbar1n}, namely, $\psi_i=\lambda_1+\delta_{0}^{\{1,2\}}$, and finally that $\lambda_1^2=0$ and $\lambda_1=\frac{1}{24} \tikz[baseline=-1mm]{\coordinate (A) at (0,0);\lp{A}{0};\legm{A}{180}{1};\gg{0}{A};}$ on $\oM_{1,1}$. All in all, we obtain
$$u_{t_1} = u u_x + \frac{\eps^2}{12}u_{xxx},$$
which we recognize as the KdV equation. It turns out (see for instance \cite[Lemma 2.4]{Bur15}) that the KdV equation determines uniquely the dispersive deformation of all other equations in the DR hierarchy upon imposing commutativity of the flows. This way the DR hierarchy construction proves for us the existence of an infinite number of symmetries for the KdV equations and its integrability.
\end{example}

\medskip

\begin{example}\label{example:ILW}
Another example of ct-F-CohFT with $V=\mbC$ is the Hodge polynomial $\Lambda(s)$. We did mention how it is compatible with the boundary stratification and this in particular implies the splitting axiom for ct-F-CohFTs, while $S_n$-equivariance is trivial. The relevant intersection numbers for the DR hierarchy this time are
$$\int_{\DR_g(-\sum_{j=1}^n a_j,0,a_1,\ldots,a_n)} \lambda_g \lambda_j \psi_2^d , \qquad 0\leq j \leq g.$$
Dimension counting this time gives us $2g-1+n=d+g+j$. Clearly, for $g=0$ the hierarchy coincides with the one for the trivial ct-F-CohFT. Let $d=1$, then for $g\geq 1$ we have $n=2-(g-j)$ and, considering that $\lambda_g^2=0$, the only case which contributes nontrivially is $n=1$ and $g-j=1$. In this case we have (see e.g. \cite{CMW12})
$$\int_{\DR_g(-a,0,a)} \lambda_g \lambda_{g-1} \psi_2=a^{2g}\frac{|B_{2g}|}{(2g)!},$$
which implies \cite{Bur15-DR-hierarchy}
$$u_{t_1}=u u_x + \sum_{g\geq 1}\eps^{2g} s^{g-1}\frac{|B_{2g}|}{(2g)!} u_{2g+1}.$$
Once again the other flows are uniquely determined by their dispersionless limit and the fact that they commute, and we recognize in the above equation the Intermediate Long Wave equation, another well known integrable system.
\end{example}

\medskip

Both of the above examples of ct-F-CohFTs are in fact F-CohFTs, by which we mean that the linear maps $c_{g,n+1}$ actually take value in $H^{\even}(\oM_{g,n})$. We shall give now an example of a strictly compact type F-CohFT.

\medskip

\begin{example}\label{example:discrete KdV}
We have seen that kappa classes behave additively upon pullback to boundary strata, i.e. $\xi_\Gamma^* \kappa_a = \sum_{v \in V(\Gamma)} p_{v}^* \kappa_a$ for any stable graph $\Gamma$. This implies that the exponential of any linear combination of kappa classes is an F-CohFT with vector space $V=\mbC$,
$$
\widetilde{c}_{g,n+1}(e^1\otimes e_1^{\otimes n})=\exp\left(\sum_{a \geq 1} s_a \kappa_a\right), \qquad s_a\in \mbC.
$$
Denote by $\mathbb{P}$ the above class $\widetilde{c}_{g,n+1}(e^1\otimes e_1^{\otimes n})$ for $s_a$ such that $\exp(-\sum_{a\geq 1}s_a t^a) = \sum_{k\geq 0} (-1)^k (2k-1)!! t^k$, and by $\mathbb{P}_d$ its part of cohomological degree $2d$. It is a result of Pixton \cite{Pix13} that $\mathbb{P}_d=0$ on $\cM_{g,n+1}^{\ct}$ for $2d>2g-1+n$ (this is not the case on $\oM_{g,n+1}$). This fact guarantees that the top degree part of $\mathbb{P}$, i.e. the class $c_{g,n+1}(e^1\otimes e_1^{\otimes n})=\mathbb{P}_{\frac{1}{2}(2g-1+n)}$ for $n$ odd and zero otherwise, is a ct-F-CohFT that is not an F-CohFT. It turns out \cite{BR21,BRR26} that the associated double ramification hierarchy is the discrete KdV hierarchy, whose primary flow is given by
$$
u_{t_0} = \d_x\left(\frac{1}{4} T(Tu)^2\right), \qquad \text{where}\quad T:=\sqrt{R}, \quad R:=\frac{2}{\mathrm{i}\eps \d_x }\frac{\Lambda_{\mathrm{i}\eps} -1}{\Lambda_{\mathrm{i}\eps} + 1},
$$
or, by the change of dependent variable $v= Tu$, $v_{t_0} = \d_x\left(\frac{1}{4} R(v^2) \right)$. 
\end{example}

\medskip

\subsection{Partial CohFTs and Hamiltonian structure}
If we look at equation \eqref{eq:DRcommct} instead of equation \eqref{eq:graphDRcomm}, we get the graphical representation 
\renewcommand{\gg}[2]{\fill[color=white] (#2) circle(4.3mm) node {\color{black}$\substack{#1}$}; \draw (#2) circle (4.3mm)}
\begin{equation}\label{eq:graphDRcomm-ham}
\sum_{\substack {g_1,g_2\geq 0\\ g_1+g_2=g\\ I\sqcup J =\{1,\ldots,n\}}} 
\sum_{\substack{ b\in \mbZ}} \ 
b\left( \tikz[baseline=-1mm]{\coordinate (A) at (0,0);\coordinate (B) at (20mm,0);\draw (B)--(A);\legm{A}{90}{0_1};\leg{A}{-60};\leg{A}{-120};\gg{\DR_{g_1}}{A};\lab{A}{-90}{6.9mm}{...};\lab{A}{-90}{9mm}{a_I};\legm{B}{90}{0_2};\leg{B}{-60};\leg{B}{-120};\gg{\DR_{g_2}}{B};\lab{B}{-90}{6.9mm}{...};\lab{B}{-90}{9mm}{a_J};\lab{B}{165}{6.8mm}{b};\lab{A}{15}{6.8mm}{-b}}\right) = 0.
\end{equation}
\renewcommand{\gg}[2]{\fill[color=white] (#2) circle(2.5mm) node {\color{black}$\substack{#1}$}; \draw (#2) circle (2.5mm)}Because we didn't single out any marking, we might get the idea that something useful would come from a more symmetric version of the notion of F-CohFT, where no insertion is special. Namely, let $V=\langle e_{\alpha}\rangle_{\alpha \in [N]}$ be a vector space endowed with a symmetric nondegenerate bilinear form $\eta=\eta^{\mu\nu}e_\mu\otimes e_\nu$, then a compact type partial CohFT is a family of linear maps $c_{g,n}\colon V^{\otimes n} \to H^\even(\cM^{\ct}_{g,n})$ that is $S_n$-equivariant and satisfies the splitting axiom
$$
\gl^* c_{g,n}(\otimes_{i\in [n]} e_{\alpha_i}) = c_{g_1,|I|+1}(\otimes_{i\in I} e_{\alpha_i}\otimes e_\mu)\eta^{\mu\nu}c_{g_2,|J|+1}(e_\nu\otimes \otimes_{j\in J} e_{\alpha_j}),
$$
with respect to $\gl=\xi_\Gamma$ with $\Gamma$ any stable graph with $n$ legs, two vertices with genus~$g_1$ and~$g_2$ and markings $I$ and $J$, respectively, $g_1+g_2=g$ and $I\sqcup J=[n]$, and one separating edge. Partial CohFTs were first introduced in \cite{LRZ15}.

\medskip

Equation \eqref{eq:graphDRcomm-ham}, when intersected with a ct-partial CohFT, becomes
\renewcommand{\gg}[2]{\fill[color=white] (#2) circle(4.3mm) node {\color{black}$\substack{#1}$}; \draw (#2) circle (4.3mm)}
\begin{equation}\label{eq:graphDRcomm-ham2}
\sum_{\substack {g_1,g_2\geq 0\\ g_1+g_2=g\\ I\sqcup J =\{1,\ldots,n\}}} 
\sum_{\substack{ b\in \mbZ}} \ 
b\left( \tikz[baseline=-1mm]{\coordinate (A) at (0,0);\coordinate (B) at (20mm,0);\draw (B)--(A);\legm{A}{90}{(0_1,\beta_1)};\leg{A}{-60};\leg{A}{-120};\gg{\DR_{g_1}}{A};\lab{A}{-90}{6.9mm}{...};\lab{A}{-90}{9mm}{(a_I,\alpha_I)};\legm{B}{90}{(0_2,\beta_2)};\leg{B}{-60};\leg{B}{-120};\gg{\DR_{g_2}}{B};\lab{B}{-90}{6.9mm}{...};\lab{B}{-90}{9mm}{(a_J,\alpha_J)};\lab{B}{165}{6.8mm}{b};\lab{A}{15}{6.8mm}{-b};\lab{A}{-15}{11.0mm}{\eta}}\right) = 0,
\end{equation}
\renewcommand{\gg}[2]{\fill[color=white] (#2) circle(2.5mm) node {\color{black}$\substack{#1}$}; \draw (#2) circle (2.5mm)}where of course the bivector $\eta$ is to be fed to the ct-partial CohFTs at each side of the edge.

\medskip

What we did is giving up on having a distinguished marking with covector insertion. It should not be surprising that the corresponding generating differential polynomial now plays the role of a function, and not a component of a vector field. Indeed, let's define
\begin{equation*}\label{eq:DR hamiltonians}
\begin{split}
&\og_{\beta,d}:=\int g_{\beta,d} dx,\\
&g_{\beta,d}:=\sum_{\substack{g,n\ge 0,\,2g+n-1>0\\k_1,\ldots,k_n\geq 0\\\sum_{j=1}^n k_j=2g}} \frac{\eps^{2g}}{n!} \Coef_{(a_1)^{k_1}\ldots(a_n)^{k_n}} \left(\int_{\DR_g(-\sum a_i,a_1,\ldots,a_n)} \hspace{-1.3cm}\lambda_g \psi_1^d c_{g,n+1} (e_\beta\otimes \otimes_{j=1}^n e_{\alpha_j}) \right)\prod_{j=1}^n u^{\alpha_j}_{k_j},
\end{split}
\end{equation*}
for $1\leq \beta\leq N$ and $d\geq 0$. We need to explain $\int \cdot \,dx$ notation: this means that $\og_{\beta,d}$ is the equivalence class of the differential polynomial $g_{\beta,d}$ in the quotient of the (degree $0$ part of) the ring $\mbC[[u^*_*,\eps]]$ by the image of the operator $\d_x$. Indeed, $\int\cdot\, dx$ is a suggestive notation for equivalence classes in this quotient because integrals with respect to $x$ are insensitive to $\d_x$-exact terms. Elements of this quotient are called local functionals. It is interesting to look at the Fourier components expression for the local functionals $\og_{\beta,d}$,
$$
\og_{\beta,d}:=\sum_{\substack{g\ge 0,\,n\ge 2\\a_1,\ldots,a_n\in \mbZ}} \frac{(-\eps^2)^{g}}{n!}  \left(\int_{\DR_g(0,a_1,\ldots,a_n)} \hspace{-1.3cm}\lambda_g \psi_1^d c_{g,n+1} (e_\beta\otimes \otimes_{j=1}^n e_{\alpha_j}) \right)\prod_{j=1}^n p^{\alpha_j}_{a_j}.
$$
We observe that integrating in $dx$ has the effect of setting the ramification profile of the first marking to $0$, so that it becomes an internal marked point, as needed for equation \eqref{eq:graphDRcomm-ham2}.

\medskip

Now, equation \eqref{eq:graphDRcomm-ham2} is exactly the commutativity equation
$$
\{\og_{\beta_1,d_1},\og_{\beta_2,d_2}\} = 0
$$
with respect to the Hamiltonian operator $\eta\d_x$, i.e. with respect to the Poisson bracket
$$\{\overline{f},\overline{h}\} := \int \left[ \frac{\delta \overline{f}}{\delta u^\mu}\eta^{\mu\nu} \d_x\left(\frac{\delta \overline{h}}{\delta u^\nu}\right)\right]dx,\qquad \text{where}\quad \frac{\delta \overline{f}}{\delta u^\alpha} = \sum_{k\geq 0} (-\d_x)^k \frac{\d f}{\d u^\alpha_k}.$$
Indeed, it's easy to compute the bracket of Fourier components of $u^*$ as $\{p^\alpha_a,p^\beta_b\}= \mathrm{i}a\eta^{\mu \nu} \delta_{a,-b}$, where we recognize the edge factor of equation \eqref{eq:graphDRcomm-ham2}.

\medskip

So we just discovered that ct-partial CohFTs produce integrable Hamiltonian evolutionary PDEs with respect to a specific Poisson bracket, while ct-F-CohFTs in general give potentially non-Hamiltonian systems. Of course ct-partial CohFTs $c_{g,n}$ give in particular ct-F-CohFts~$\widetilde{c}_{g,n+1}$ as $\widetilde{c}_{g,n+1}(e^\alpha\otimes \otimes_{i\in [n]} e_{\alpha_i}):=\eta^{\alpha\mu}c_{g,n+1}(e_\mu\otimes \otimes_{i \in [n]}e_{\alpha_i})$. At the level of generating functions, this gives exactly the relation between flows and Hamiltonians
$$u^\alpha_{t^\beta_d} = \eta^{\alpha\mu}\d_{x}\left(\frac{\delta \overline{g}_{\beta,d}}{\delta u^\mu}\right).$$

\medskip

\begin{example}
The F-CohFTs of Examples \ref{example:KdV}, \ref{example:ILW}, and \ref{example:discrete KdV} actually come from partial CohFTs $c_{g,n}(e_1^{\otimes n}) = 1$, $c_{g,n}(e_1^{\otimes n}) = \Lambda(s)$ and $\widetilde{c}_{g,n}(e_1^{\otimes n})= \exp\left(\sum_{a \geq 1} s_a \kappa_a\right)$, respectively, all of them with metric $\eta_{11} = 1$, so all of these systems are actually Hamiltonian with respect to the Hamiltonian operator $\d_x$. For instance, for $c_{g,n}(e_1^{\otimes n}) = 1$ we get
$$\overline{g}_{1}=\int\left[\frac{u^3}{6} + \frac{\eps^2}{24} u u_{xx} \right] dx.$$
\end{example}

\medskip

\begin{example}
Let us give an example of a rank $1$ F-CohFT that is not a partial CohFT. The idea is starting from the family of partial CohFTs  $\widetilde{c}_{g,n+1}(e_1^{\otimes n})= \exp\left(\sum_{a \geq 1} s_a \kappa_a\right)$ on $\cM^{\ct}_{g,n+1}$ and break its $S_{n+1}$-symmetry by adding psi classes at the marked points with coefficients which single out the first marking as special. As a first attempt we could write something like
$$\exp\left(\sum_{a \geq 1} s_a \left(\kappa_a + (-1)^a \psi_1^a - \sum_{i=2}^{n+1}\psi_i^a\right)\right)$$
and, for even $a$, the first marked point would be treated differently from the others, because of the sign difference. Of course this arbitrary modification ruins the splitting axiom, as the extra psi terms will fail to appear at the separating node, but we can salvage it by adding the corresponding counterterm, namely 
$$
c_{g,n+1}(e_1^{\otimes n})= \exp\Bigg(\sum_{a \geq 1} s_a \Bigg(\kappa_a+ (-1)^a \psi_1^a - \sum_{i=2}^{n+1}\psi_i^a + \sum_{\substack{i+j=a-1\\ g_1+g_2=g \\ I\sqcup J = \{2,\ldots,n+1\}}} (-1)^j
\tikz[baseline=-1mm]{\coordinate (A) at (0,0);\coordinate (B) at (16mm,0);\draw (A)--(B);\legm{A}{180}{1};\leg{A}{-60};\leg{A}{-120};\leg{B}{-60};\leg{B}{-120};\gg{g_1}{A};\gg{g_2}{B};\lab{A}{-90}{6.9mm}{...};\lab{A}{-90}{10mm}{\underbrace{\phantom{aaa}}_{\text{$I$}}};\lab{B}{-90}{6.9mm}{...};\lab{B}{-90}{10mm}{\underbrace{\phantom{aaa}}_{\text{$J$}}};\lab{A}{25}{4.8mm}{\psi^i};\lab{B}{155}{4.8mm}{\psi^j}}
\Bigg)\Bigg).
$$
The edge term we added can equivalently be written as $\sum_{\Gamma} (\xi_\Gamma)_*\left(\frac{\psi_{e'}^a-(-\psi_{e''})^a}{\psi_{e'} + \psi_{e''}}\right)$, where the sum runs over all stable graphs with two vertices and one separating edge $e=(e',e'')$ such that the first marking is on the vertex to which the half edge $e'$ is attached. Recalling that the first Chern class of the normal bundle to $\xi_\Gamma$ is $-(\psi_{e'}+\psi_{e''})$ we see that, upon pulling back along~$\xi_\Gamma$, the extra term produces $-(\psi_{e'}+\psi_{e''})\left(\frac{\psi_{e'}^d-(-\psi_{e''})^a}{\psi_{e'} + \psi_{e''}}\right)$ which is the correct counterterm needed for the splitting axiom. Comparing with formula \eqref{eq:Mumford} we have a geometric interpretation for the summands with odd $a$ in the exponential, for which the $S_{n+1}$-symmetry is preserved, as $s_a \frac{(a+1)!}{B_{a+1}} \mathrm{ch}_{a}(\mathbb{E})$ (the factor of $\frac{1}{2}$ in front of the sum in formula \eqref{eq:Mumford} has been replaced by the condition that the first marking is on the vertex to which $e'$ is attached, while $|\Aut(\Gamma)|=1$ for these stable graphs). On the other hand, the terms with even $a$, beside being clearly only $S_n$-symmetric in the last $n$ markings, present an edge term which changes sign upon exchanging~$e'$ with~$e''$.
\end{example}

\medskip

\subsection{CohFTs and quantization}
Equation \eqref{eq:DRcommct} on $\cM^{\ct}_{g,n}$ that we used to deduce commutativity of the DR hierarchy's Hamiltonians in the previous section, stems from a more general equation, namely equation \eqref{eq:DRcomm} on $\oM_{g,n}$. We did not use that until now because we were looking for equations corresponding to the vanishing of Lie or Poisson brackets, and these are bidifferential operators of order at most $1$ in the arguments, which, in our graphical representation, correspond to a single edge connecting two vertices. In order to go from equation \eqref{eq:DRcomm} to equation \eqref{eq:DRcommct} we had to multiply by $\lambda_g$, which excluded multiple edges between the two vertices.

\medskip

However, if we accept to replace Poisson brackets with quantum commutators, we can use the more general equation \eqref{eq:DRcomm} to obtain a quantum integrable Hamiltonian system instead! These quantum brackets are expressed as commutators
$[\overline{f},\overline{h}] = \frac{1}{\hbar}(\overline{f}\star\overline{h} - \overline{h}\star\overline{f})$
with respect to a star product whose expression in terms of the Fourier components of $u*$ is
$$ 
\overline{f} \star \overline{g} =\overline{f} \exp\left(\sum_{k\geq 0}\mathrm{i}\hbar b \overleftarrow{\frac{\d}{\d p^\mu_b}} \eta^{\mu\nu} \overrightarrow{\frac{\d}{\d p^\nu_{-b}}}\right)\overline{g}.
$$
The above formula is the Moyal quantization of the DR hierarchy's Poisson bracket and we see bidifferential operators of any order appearing, with the correct edge factor to match the one in formula \eqref{eq:DRcommct}.

\medskip

Since equation \eqref{eq:DRcomm} is valid on $\oM_{g,n}$ we need our linear maps $c_{g,n}$ to take vale in $H^*(\oM_{g,n})$, and since equation \eqref{eq:DRcomm} contains gluing maps with respect to stable graphs with two vertices connected by multiple edges, we need the $c_{g,n}$ to satisfy the sequence of splitting axioms,
\begin{equation}\label{eq:bananasplit}
\gl^* c_{g,n}(\otimes_{i\in [n]} e_{\alpha_i}) = c_{g_1,|I|+p}(\otimes_{i\in I} e_{\alpha_i}\otimes \otimes_{r=1}^p e_{\mu_r})\left(\prod_{r=1}^p\eta^{\mu_r\nu_r}\right)c_{g_2,|J|+p}(\otimes_{r=1}^p e_{\nu_r}\otimes \otimes_{j\in J} e_{\alpha_j}),
\end{equation}
with respect to $\gl=\xi_\Gamma$ with $\Gamma$ any stable graph with $n$ legs, two vertices with genus $g_1$ and $g_2$ and markings $I$ and $J$, respectively, $p$ edges connecting the two vertices, $g_1+g_2+p-1=g$ and $I\sqcup J=[n]$, for any integer $p\geq 1$.

\medskip

A special case of such systems of linear maps is given by cohomological field theories (CohFTs). These are partial CohFTs (taking values in $H^{\even}(\oM_{g,n})$) satisfiying just one extra gluing axiom (sometimes called the loop axiom):
$$
\gl^* c_{g+1,n}(\otimes_{i\in [n]} e_{\alpha_i}) = c_{g,n+2}(\otimes_{i\in [n]} e_{\alpha_i}\otimes  e_{\mu} \otimes e_\nu)\eta^{\mu\nu},
$$
with respect to $\gl=\xi_\Gamma$ with $\Gamma$ the stable graph with one vertex, $n$ legs, and one loop (self-edge). Clearly, combining this loop gluing axiom with the splitting axiom at separating edges one can deduce the sequence of axioms \eqref{eq:bananasplit}. CohFTs were first introduced in \cite{KM94}.

\medskip

In the presence of either a CohFT or a partial CohFT satisfying the extra axioms \eqref{eq:bananasplit}, we can define quantum DR Hamiltonians \cite{BR16b},
$$
\overline{G}_{\beta,d}:=\sum_{\substack{g,n\geq 0,\,2g+n-1>0\\a_1,\ldots,a_n\in \mbZ}} \frac{(\mathrm{i\hbar})^{g}}{n!}  \left(\int_{\DR_g(0,a_1,\ldots,a_n)} \Lambda\left(\frac{-\eps^2}{\mathrm{i}\hbar}\right)  \psi_1^d c_{g,n+1} (e_\beta\otimes \otimes_{j=1}^n e_{\alpha_j}) \right)\prod_{j=1}^n p^{\alpha_j}_{a_j},
$$
for $1\leq \beta\leq N$ and $d\geq 1$.

\medskip

Notice that the role of the class $\Lambda(\frac{-\eps^2}{\mathrm{i}\hbar})$ is merely to provide the correct classical limit, as $\hbar\to 0$, to the classical DR hamiltonians, but we could have omitted that factor altogether (after all $\Lambda(s)$ is a CohFT and is invertible as a cohomology class, and CohFTs form a ring under cup product in cohomology and direct sum of vector spaces). By equation $\eqref{eq:DRcomm}$, these commute with respect to the aforementioned quantum brackets.

\medskip

Notice also that, for computing the quantum DR hierarchy, Hain's formula \eqref{eq:Hain} is not enough anymore: we need Pixton's full formula since $\lambda_g$ is not there anymore and the non compact type components of $\DR_g(0,A)$ now matter.

\medskip

\begin{example}
The partial CohFTs $c_{g,n}(e_1^{\otimes n}) = 1$, $c_{g,n}(e_1^{\otimes n}) = \Lambda(s)$, and $\widetilde{c}_{g,n}(e_1^{\otimes n})= \exp\left(\sum_{a \geq 1} s_a \kappa_a\right)$ of Examples \ref{example:KdV}, \ref{example:ILW}, and the beginning of \ref{example:discrete KdV}, respectively, are actually CohFTs, as they all satisfy the loop axiom too, so all of these systems possess a quantization. For instance, for $c_{g,n}(e_1^{\otimes n}) = 1$ we get
$$
\overline{G}_{1}=\int\left[\frac{u^3}{6} + \frac{\eps^2}{24} u u_{xx}  - \frac{\mathrm{i\hbar}}{24} u \right]dx.
$$
The extra term, linear in $\hbar$, comes from $g=1$ and $n=1$, namely,
$$
\int_{\DR_1(0,0)}\psi_1 = \frac{1}{24},
$$
where we have used that $\DR_g(0,\ldots,0)=(-1)^g\lambda_g$. Notice that, instead, quantizability is not a priori given for the class $c_{g,n}( e_{1}^{\otimes n})= \mathrm{P}_{\frac{1}{2}(2g-1+n)}$ of Example~\ref{example:discrete KdV}.
\end{example}

\medskip

\subsection{Flat units}

A final aspect we can consider does not directly influence the possibility of constructing classical or quantum commuting systems and comes more from imposing a natural compatibility of our families of classes $c_{g,n}$ with the forgetful morphisms $\pi$, the maps between moduli spaces which forget one of the marked points.

\medskip

Compact type or not, we say that F-CohFTs, partial CohFTs or CohFTs have a (flat) unit if they satisfy one extra axiom: the vector space $V$ comes equipped with a special vector $e\in V$, called the unit, such that, for F-CohFTs
$$
\pi^* c_{g,n+1}(e^\alpha \otimes\otimes_{i\in [n]} e_{\alpha_i})= c_{g,n+2}(e^\alpha \otimes \otimes_{i\in [n]} e_{\alpha_i}\otimes e)
$$
for $g,n\geq 0$ and $2g-1+n>0$, and
$$
\delta^\mu_\nu = c_{0,3}(e^\mu\otimes e_\nu \otimes e),
$$
while for partial CohFTs and CohFTs
$$\pi^* c_{g,n}(\otimes_{i\in [n]} e_{\alpha_i})= c_{g,n+1}(\otimes_{i\in [n]} e_{\alpha_i}\otimes e)$$
for $g,n\geq 0$ and $2g-2+n>0$, and 
$$\eta_{\mu\nu}=c_{0,3}(e_\mu\otimes e_{\nu}\otimes e),$$
with $(\eta_{\mu \nu})$ the inverse matrix to $(\eta^{\mu\nu})$.\\

\medskip

The main effect of the unit is that the following equation is valid for the flows of the DR hierarchy
$$\frac{\d P^\alpha_{\beta,d}}{\d u^e} = P^\alpha_{\beta,d-1}, \qquad d \geq 1,$$
where, if  $e=A^\alpha e_\alpha$, $A^\alpha \in \mbC$, then $\frac{\d}{\d u^e} := A^\alpha \frac{\d}{\d u^\alpha}$ and
$$P^\alpha_{e,0} =  u^\alpha, $$
where $P^\alpha_{e,0}:=A^\mu P^\alpha_{\mu,0}$.
Similarly, in the Hamiltonian case,
$$\frac{\d g_{\beta,d}}{\d u^e} = g_{\beta,d-1}, \qquad d\geq 1,$$
and
$$\overline{g}_{e,0} = \int \frac{1}{2} \eta_{\mu\nu}u^\mu u^\nu \,dx,$$
where $\overline{g}_{e,0}:=A^\mu \overline{g}_{\mu,0}$.
These come directly from the fact that for any differential polynomial $f$, $\frac{\d f}{\d u^\alpha} = \frac{\d f}{\d p^\alpha_0}$ and, in the integrals over $\oM_{g,n}$ making up the generating functions $\frac{\d P^\alpha_{\beta,d}}{\d u^e}$ or $\frac{\d g_{\beta,d}}{\d u^e}$, every integrand but the psi class is a pullback with respect to the morphism $\pi$ forgetting one of the markings (the one corresponding to the $p^\alpha_0$ derivative). Applying the push-pull formula we end up lowring the power of the psi class by one.

\smallskip
 
Another consequence of the presence of a flat unit is that there is no difference between axiom \eqref{eq:bananasplit} and the loop axiom of a CohFT as one can see by composing the pullback with respect to the forgetful and gluing morphisms.

\medskip

\begin{example}
While the F-CohFTs (actually coming from CohFTs) of Examples \ref{example:KdV} and \ref{example:ILW} have flat unit $e = e_1$, the F-CohFTs (actually coming from partial CohFTs) mentioned at the beginning of Example \ref{example:discrete KdV} do not have a flat unit.
\end{example}

\end{document}